\UseRawInputEncoding
\documentclass[11pt, reqno]{amsart}
\usepackage{amsfonts,latexsym,enumerate}
\usepackage{amsmath}
\usepackage{amscd}

\usepackage{float,amsmath,amssymb,mathrsfs,bm,multirow,graphics}
\usepackage[dvips]{graphicx}
\usepackage[percent]{overpic}
\usepackage[numbers,sort&compress]{natbib}
\usepackage{todonotes}
\usepackage{grffile}

\usepackage{tikz}

\newcommand{\R}{{\Bbb R}}

\newcommand{\C}{{\Bbb C}}

\newcommand{\diag}{\text{\upshape diag\,}}

\newcommand{\im}{\text{\upshape Im\,}}

\newcommand{\sol}{\text{\upshape sol}}

\newcommand{\lead}{\text{\upshape lead}}

\newcommand{\noshift}{\text{\upshape no shift}}

\DeclareMathOperator{\sech}{sech}

\DeclareMathOperator{\res}{Res}

\def\XXint#1#2#3{{\setbox0=\hbox{$#1{#2#3}{\int}$}
\vcenter{\hbox{$#2#3$}}\kern-.5\wd0}}

\newtheorem{theorem}{Theorem}[section]

\newtheorem{remark}[theorem]{Remark}

\newtheorem{figuretext}[theorem]{Figure}
\newtheorem{RHproblem}[theorem]{RH problem}

\numberwithin{equation}{section}

\usepackage[colorlinks=true]{hyperref}
\hypersetup{urlcolor=blue, citecolor=red, linkcolor=blue}

\input epsf
\title[Numerical scheme for the ``bad'' Boussinesq equation]
{Numerical scheme for the solution of the ``bad'' Boussinesq equation}

\author{Christophe Charlier}
\address{CC: Centre for Mathematical Sciences, Lund University, 221 00 Lund, Sweden.}
\email{christophe.charlier@math.lu.se}

\author{Daniel Eriksson}
\address{DE: Department of Mathematics, KTH Royal Institute of Technology, 100 44 Stockholm, Sweden.}
\email{deri3@kth.se}

\author{Jonatan Lenells}
\address{JL: Department of Mathematics, KTH Royal Institute of Technology, 100 44 Stockholm, Sweden.}
\email{jlenells@kth.se}

\begin{document}

\begin{abstract}
We present a numerical scheme for the solution of the initial-value problem for the ``bad'' Boussinesq equation. 
The accuracy of the scheme is tested by comparison with exact soliton solutions as well as 
with recently obtained asymptotic formulas for the solution.
%Should we also mention that our scheme gives a significantly smaller $L^{\infty}$-error than previous numerical works on the one-soliton?
\end{abstract}

\maketitle



%\setcounter{tocdepth}{1} 
%\tableofcontents

\section{Introduction}
In this paper, we present a numerical scheme for the solution of the initial-value problem for the ``bad'' Boussinesq equation
\begin{align}\label{badboussinesq}
u_{tt} - u_{xx} - (u^2)_{xx} - u_{xxxx} = 0,
\end{align}
where subscripts denote partial derivatives and $u(x,t)$ is a real-valued function of a space variable $x$ and a time variable $t$.
Equation (\ref{badboussinesq}) was derived a century and a half ago as a model for the propagation of long water waves of small amplitude by J. Boussinesq \cite{B1872}. Along with the Korteweg--de Vries (KdV) equation, (\ref{badboussinesq}) is one of the fundamental models for dispersive water waves in the weakly nonlinear regime \cite{J1997}. Just like the KdV equation, equation (\ref{badboussinesq}) admits a Lax pair formulation \cite{Z1974} and supports soliton solutions \cite{H1973}. However, in contrast to the unidirectional KdV equation, (\ref{badboussinesq}) is a bidirectional equation, i.e., it models waves traveling in both the right and left directions.

Equation (\ref{badboussinesq}) is also known as the ``nonlinear string equation'' \cite{Z1974}. In addition to water waves and nonlinear strings, it also models other physical phenomena, including the propagation of ion-sound waves in a plasma \cite{M1978}, the dynamics of the anharmonic lattice in the Fermi--Pasta--Ulam problem \cite{M1978}, and quite generally nonlinear lattice waves in the continuum limit \cite{T1975}. 

A distinguishing feature of equation (\ref{badboussinesq}) is that it is ill-posed. At the linear level, the ill-posedness is easily observed by substituting the plane wave $u(x,t) = e^{i(kx - \omega t)}$, $k \in \R$, into the linearized version of (\ref{badboussinesq}). This gives the dispersion relation $\omega = \pm \sqrt{k^2(1 - k^2)}$. Since $\omega$ has a nonzero imaginary part whenever $|k| > 1$, Fourier modes with $|k| > 1$ grow (or decay) exponentially in time. This feature makes the initial value problem linearly ill-posed. A deeper analysis shows that the ill-posedness persists at the nonlinear level: in fact, it was shown in \cite{CLmain} that there exist Schwartz class solutions of (\ref{badboussinesq}) that blow up at time $T$ for any $T \in (0, +\infty)$. Moreover, a wide range of asymptotic scenarios can occur as $T$ is approached \cite{Cbad}.

Due to its ill-posedness, equation (\ref{badboussinesq}) is notoriously difficult to study both analytically and numerically (see e.g. \cite{DH1999, D2008}). Accordingly, the terminology ``bad'' Boussinesq equation was introduced in \cite{KL1977} for (\ref{badboussinesq}). This is in contrast to the so-called ``good'' Boussinesq equation, which has the same form as (\ref{badboussinesq}) except that the term $u_{xxxx}$ has the opposite sign. Changing the sign of the $u_{xxxx}$-term makes the equation linearly well-posed, and a large number of works are dedicated to the study of the ``good'' Boussinesq equation: for example, numerical investigations include \cite{MMM1985, OS1990, PS1997, E2003, IB2003, BTN2005, B1998, B2007, IM2014}, asymptotic estimates and spectral results are obtained in \cite{LS1995, F2011, L1997, CL2022, CLW2022}, local existence and well-posedness results can be found in \cite{BFH2019, BS1988, F2009, F2011, HM2015, KT2010, TM1991, X2006, CT2017}, and global existence results are derived in \cite{BS1988, L1993, TM1991, X2006}. The ``good'' Boussinesq equation also has various applications, but it is the version derived by J. Boussinesq (i.e., the ``bad'' version) that is relevant for water waves. Since it is relatively straightforward to solve the initial value problem for the ``good'' Boussinesq equation numerically to high accuracy \cite{MMM1985}, we henceforth focus on the ``bad'' Boussinesq equation (\ref{badboussinesq}).

The purpose of this paper is to propose a numerical scheme for the solution of the initial-value problem for (\ref{badboussinesq}). The proposed scheme involves six steps. The first step is to rewrite the equation (\ref{badboussinesq}), which is second-order in $t$, as a first-order system. The second step is to apply a discrete Fourier transform to this first-order system, and thus obtain a system of ordinary differential equations that describes the time evolution of the Fourier modes. The third step is to set all Fourier modes of frequency above a carefully chosen cut-off frequency to zero. The fourth step is to add a damping term to one of the equations in the first-order system so that the transition from the nonzero Fourier modes to the vanishing ones is smoothened out. The fifth step is to solve the damped system of ordinary differential equations by means of, for example, a Runge--Kutta method, to obtain approximate Fourier coefficients. In the sixth and last step, an approximate solution of (\ref{badboussinesq}) is obtained by applying an inverse Fourier transform.

Our construction of the above numerical scheme was inspired by recent developments in the analytical understanding of the ``bad'' Boussinesq equation. More precisely, in \cite{CLmain}, the flow of (\ref{badboussinesq}) was linearized by means of a ``nonlinear Fourier transform''. This led to existence, uniqueness, and blow-up results, but also to the identification of a physically meaningful class of global solutions. The global solutions arise as follows. As already mentioned, the linear ill-posedness of (\ref{badboussinesq}) leads to the exponential amplification of high frequency Fourier modes. 
%In particular, solutions will grow exponentially unless the unstable modes vanish identically. 
The presence of unstable modes does not contradict Boussinesq's derivation of (\ref{badboussinesq}) as a model for water waves, because the derivation assumes that the waves are long, in the sense that the wavelength is long compared to the depth of the water (see Section \ref{regimesec} for more details). Thus, physically relevant solutions arise from initial data without high-frequency modes. However, since the equation is nonlinear, high-frequency modes will be generated by the time evolution, even if none are present initially. The ``nonlinear Fourier transform'' of \cite{CLmain} makes it possible to identify the unstable modes in a nonlinearly precise way. In particular, it is shown in \cite{CLmain} that initial data with vanishing ``nonlinear high-frequency modes'', give rise to global solutions, which can be expected to describe the physical characteristics of water waves according to Boussinesq's model. The numerical scheme proposed here grew out of an attempt to capture these physically relevant solutions of (\ref{badboussinesq}). In particular, the scheme time-evolves the Fourier coefficients of the solution while setting Fourier modes of high enough frequency to zero in a smooth way.

%The choice of the cut-off frequency is important; we show that it can be chosen so that (i) numerical instabilities originating from the linear ill-posedness are avoided and (ii) solutions of (\ref{badboussinesq}) that describe small-amplitude long waves are well approximated. Since equation (\ref{badboussinesq}) is derived under the assumption that waves are long and of small amplitude 
Being integrable, equation (\ref{badboussinesq}) admits large families of exact soliton solutions \cite{H1973}. We test our numerical scheme by comparing the approximate solution it generates with the exact solution in the case of one-solitons. We also test the scheme by comparing the approximate solution it generates with the asymptotic formulas presented in \cite{CLmain}. These asymptotic formulas describe the large time behavior of the solution of (\ref{badboussinesq}) and allow us to test the scheme for more general initial data.

\subsection{Relation to earlier work}
We are aware of six earlier papers dealing with numerical schemes for (\ref{badboussinesq}), see \cite{DH1999, IB2003, BTN2005, B1998, B2007, UEK2021}.
In \cite{B1998}, a numerical solution of (\ref{badboussinesq}) was obtained by using the so-called method of lines. A finite difference scheme for the numerical solution of (\ref{badboussinesq}) was proposed and investigated in \cite{DH1999}, where accuracy was improved by the application of filtering and regularization techniques. 
In \cite{IB2003}, a fourth-order in time and second-order in space finite-difference scheme was developed, and a predictor-corrector pair was proposed for the solution of the resulting non-linear system.
The work \cite{BTN2005} introduces linearized numerical schemes to study (\ref{badboussinesq}), and a predictor-corrector scheme was presented in \cite{B2007}. In \cite{UEK2021}, a Galerkin finite element method is applied to (\ref{badboussinesq}). Common to the schemes in \cite{DH1999, IB2003, BTN2005, B1998, B2007}, is that they time-evolve the solution in physical space. The scheme proposed in this paper, in contrast, time-evolves the solution in spectral (Fourier) space. 

The schemes proposed in \cite{IB2003, BTN2005, B1998, B2007, UEK2021} are all tested on the one-soliton of amplitude $0.369$ and numerical $L^\infty$-errors are given. 
The amplitude $0.369$ is too large to make the solution comfortably fit within the small-amplitude long-wave regime where the Boussinesq equation is expected to be a good model for water waves, and where our numerical scheme is expected to provide a good approximation. Nevertheless, in Section \ref{solitonsec}, we have applied our scheme to the soliton of amplitude $0.369$ and compared the resulting error with the errors found in \cite{IB2003, BTN2005, B1998, B2007, UEK2021}. 
The $L^\infty$-errors we find at $t = 72$ are roughly five times smaller than those found in \cite{B1998, BTN2005, B2007}. In \cite{IB2003}, the $L^\infty$-error $0.074056$ is reported at $t = 50$, whereas we find $0.00734$. In \cite[Table 9]{UEK2021}, the $L^\infty$-error at $t = 72$ is reported to be $0.0207$, whereas we get $0.00880$. For the soliton of amplitude $0.1$, the $L^\infty$-errors $0.262 \times 10^{-2}$ and $0.147 \times 10^{-3}$ at $t = 36$ are reported in \cite[Table 3]{B1998} and \cite[Table 9]{UEK2021}, respectively, while we find $0.212 \times 10^{-4}$. In general, our method yields significantly smaller $L^\infty$-errors than those reported in \cite{IB2003, BTN2005, B1998, B2007, UEK2021} except at very small times.
See Section \ref{solitonsec} for further details.

%\cite[Figure 2]{B1998} has error in sup norm about 0.04 for t=72, we get just above 0.01
%about 0.04-0.05 also in \cite[Figure 1]{BTN2005}
%about 0.04 also in \cite[Figure 1b]{B2007}.
%
%\cite[Table 2]{IB2003} only lists times up until $t = 50$. At $t = 50$, they already have $e \approx 0.074$.

The numerical scheme is presented in Section \ref{schemesec}. In Section \ref{regimesec}, we comment on the physical regime where (\ref{badboussinesq}) is expected to model water waves. In Sections \ref{solitonsec} and \ref{asymptoticsec}, the scheme is tested by comparison with exact soliton solutions and with the asymptotic formulas of \cite{CLmain}, respectively.

\section{A numerical scheme}\label{schemesec}
In this section, we describe the numerical scheme proposed in this paper. We consider spatially periodic solutions of (\ref{badboussinesq}) of period $2L$, i.e., solutions $u(x,t)$ such that $u(x + 2L,t) = u(x,t)$, for some fixed $L>0$. 
We assume that $u_0(x) := u(x,0)$ and $u_1(x) := u_t(x,0)$ are given, and that 
\begin{align}\label{u1zeromean}
\int_{-L}^L u_1(x) dx = 0.
\end{align}
 The condition (\ref{u1zeromean}) is fulfilled for any physically meaningful data for which the total mass $\int u dx$ is conserved, because (\ref{badboussinesq}) implies that 
\begin{align*}
\frac{d^2}{dt^2}\int_{-L}^L u dx = 0, \quad \text{i.e.} \quad \int_{-L}^L u dx = \bigg(\int_{-L}^L u_1 dx\bigg)t + \int_{-L}^L u_0 dx,
\end{align*}
where the boundary terms vanish because of periodicity. 

Equation (\ref{badboussinesq}) can be equivalently written as the system
\begin{align}\label{boussinesqsystem}
& \begin{cases}
 u_t = v_x,
 	\\
 v_{t} = u_{x} + (u^2)_{x} + u_{xxx},
\end{cases}
\end{align}
where, thanks to (\ref{u1zeromean}), $v(x,t) := \int_{-L}^x  u_t(x',t) dx'$ is periodic of period $2L$.

\begin{remark}
Our scheme is formulated for spatially periodic solutions of period $2L$. To treat the non-periodic initial value problem on the line, we choose $L$ so large that the initial data is appreciably different from zero only in an interval contained in $[-L, L]$. We then use the periodic extension of the restriction of the initial data to $[-L,L]$ as initial data for the periodic problem. Finally, we run the simulation only for so long that $u(x,t)$ remains approximately zero for $x$ near the endpoints $\pm L$. By choosing $L$ large enough, this leads to an accurate numerical solution also for the problem on the line.
\end{remark}

\subsection{A Fourier transformed system}
Let $N \geq 1$ be an integer and define $x_j := jL/N$, so that $\{x_j \,|\,  j = -N, -N+1, \dots, N-1\}$ is a partition of the interval $[-L,L]$ with mesh size $L/N$. For a function $f$ defined at each of the points $x_j$, we define the discrete Fourier transform $\hat{f}(j)$ for $j = -N, \dots, N-1$ by
$$\hat{f}(j) := \frac{1}{2N} \sum_{l = -N}^{N-1} f(x_{l}) e^{-\frac{\pi i j}{L} x_{l}}.$$

The inverse transform
$$F(x) := \sum_{j=-N}^{N-1} \hat{f}(j) e^{\frac{\pi i j}{L} x}$$
is an approximation of $f$ such that $F(x_j) = f(x_j)$ for each $j$. The $n$th derivative $F^{(n)}(x)$ is an approximation of $f^{(n)}(x)$ given by
$$F^{(n)}(x) = \sum_{j=-N}^{N-1} \widehat{F^{(n)}}(j) e^{\frac{\pi i j}{L} x}$$
with Fourier coefficients $\widehat{F^{(n)}}(j) := (\frac{\pi i j}{L})^n \hat{f}(j)$.
An application of a discrete Fourier transform therefore suggests that we approximate (\ref{boussinesqsystem}) by the following system of $4N$ ordinary differential equations:
\begin{align}\label{boussinesqsystemFourier}
& \begin{cases}
 \hat{U}_t(j,t) = \frac{\pi i j}{L}\hat{V}(j,t),
	\\
 \hat{V}_{t}(j,t) = \frac{\pi i j}{L}\hat{U}(j,t) + 2(\widehat{U_{x}} * \hat{U})(j,t) + (\frac{\pi i j}{L})^3 \hat{U}(j,t),
\end{cases} \quad j = -N, \dots, N-1,
\end{align}
where the convolution $\widehat{U_{x}} * \hat{U}$ is defined by
$$(\widehat{U_{x}} * \hat{U})(j,t)
= \begin{cases} 
\sum_{l=-N}^{N+j} \widehat{U_{x}}(l,t) \hat{U}(j-l,t) & \text{for $j = -N, \dots, -1$}, \\
\sum_{l=j +1-N}^{N-1} \widehat{U_{x}}(j-l,t) \hat{U}(l,t) & \text{for $j = 0, \dots, N-1$}, 
\end{cases}
$$
with $\widehat{U_{x}}(l,t) := \frac{\pi i j}{L}\hat{U}(l,t)$.
Solving (\ref{boussinesqsystemFourier}) with initial conditions given by the discrete Fourier transforms of $u(\cdot,0)$ and $v(\cdot,0)$, i.e.,
\begin{align}\label{hatUVinitialconditions}
\hat{U}(j,0) = \frac{1}{2N}\sum_{l = -N}^{N-1} u(x_{l}, 0) e^{-\frac{\pi ij}{L} x_{l}}, \qquad
\hat{V}(j,0) = \frac{1}{2N}\sum_{l = -N}^{N-1} v(x_{l}, 0) e^{-\frac{\pi ij}{L} x_{l}},
\end{align}
we expect the inverse Fourier transform
\begin{align}\label{Udef}
U(x,t) := \sum_{j=-N}^{N-1} \hat{U}(j,t) e^{\frac{\pi i j}{L} x}
\end{align}
to be a good approximation of the solution $u(x,t)$ of (\ref{badboussinesq}). 
In the case of the ``good'' Boussinesq equation, a scheme of this form indeed gives an accurate solution, and the accuracy can be increased by choosing $N$ large. However, in the case of the ``bad'' Boussinesq equation (\ref{badboussinesq}), the above scheme is only accurate if $N$ is appropriately chosen: if $N$ is too large the linear ill-posedness makes the scheme exceedingly unstable, while if $N$ is too small, the approximation fails because it ignores relevant frequencies.

\subsection{The choice of $N$}
To understand how to choose $N$ for (\ref{badboussinesq}), we observe that Fourier transformation of the linearized version of (\ref{badboussinesq}), $u_{tt} = u_{xx} + u_{xxxx}$, leads to
$$\hat{U}_{tt}(j, t) = \bigg(- \Big(\frac{\pi j}{L} \Big)^2 + \Big(\frac{\pi j}{L}\Big)^4\bigg) \hat{U}(j, t),$$
and hence
$$\hat{U}(j, t) = c_1(j)e^{ i t \sqrt{(\frac{\pi j}{L})^2 - (\frac{\pi j}{L})^4}} + c_2(j)e^{- i t \sqrt{(\frac{\pi j}{L})^2 - (\frac{\pi j}{L})^4}}$$
for some time-independent functions $c_1$ and $c_2$. Thus exponential growth of the coefficient $\hat{U}(j, t)$ is avoided provided that
\begin{align}\label{boundedcondition}
\Big(\frac{\pi j}{L}\Big)^2 - \Big(\frac{\pi j}{L}\Big)^4 \geq 0.
\end{align}
This suggests that we can achieve good numerical accuracy if $\hat{U}(j, t)$ is nonzero only when $|j| \leq \frac{L}{\pi}$. We therefore choose $N$ to be the integer part of $L/\pi$:
\begin{align}\label{Ndef}
N = \bigg\lfloor \frac{L}{\pi} \bigg\rfloor,
\end{align}
meaning that high-frequency modes with $|j| > \frac{L}{\pi}$ are set to zero.
In particular, no frequency mode with wave length below $2\pi$ can be captured by our scheme. If such high-frequency modes are present in the initial data, they will be immediately lost by our scheme. However, equation (\ref{badboussinesq}) is a good model for water waves only in the long wave regime, and physically meaningful solutions arise when no high-frequency modes are present, see Section \ref{regimesec}. More precisely, it is the ``nonlinear high-frequency modes'' mentioned in the introduction that should be absent, and (\ref{boundedcondition}) is a natural linear approximation of this nonlinear condition.

\subsection{Smooth damping out of high-frequency modes}\label{dampingsubsec}
The next step of our numerical scheme is to smoothen out the transition from the nonzero Fourier modes with $j \in [-N, N-1]$ to the vanishing modes with $j \notin [-N, N-1]$. 
It turns out that the accuracy of the scheme can be improved substantially if instead of leaving all the modes with $j \in [-N, N-1]$ unchanged and setting the modes with $j \notin [-N, N-1]$ to zero, one introduces a gradually increasing damping of the Fourier modes as $j$ approaches the endpoints $-N$ and $N-1$. 

%j_0 = 1 when j = -N, so j_0 = j + N + 1

\begin{figure}
\bigskip\begin{center}
\hspace{-.4cm}
\begin{overpic}[width=.4\textwidth]{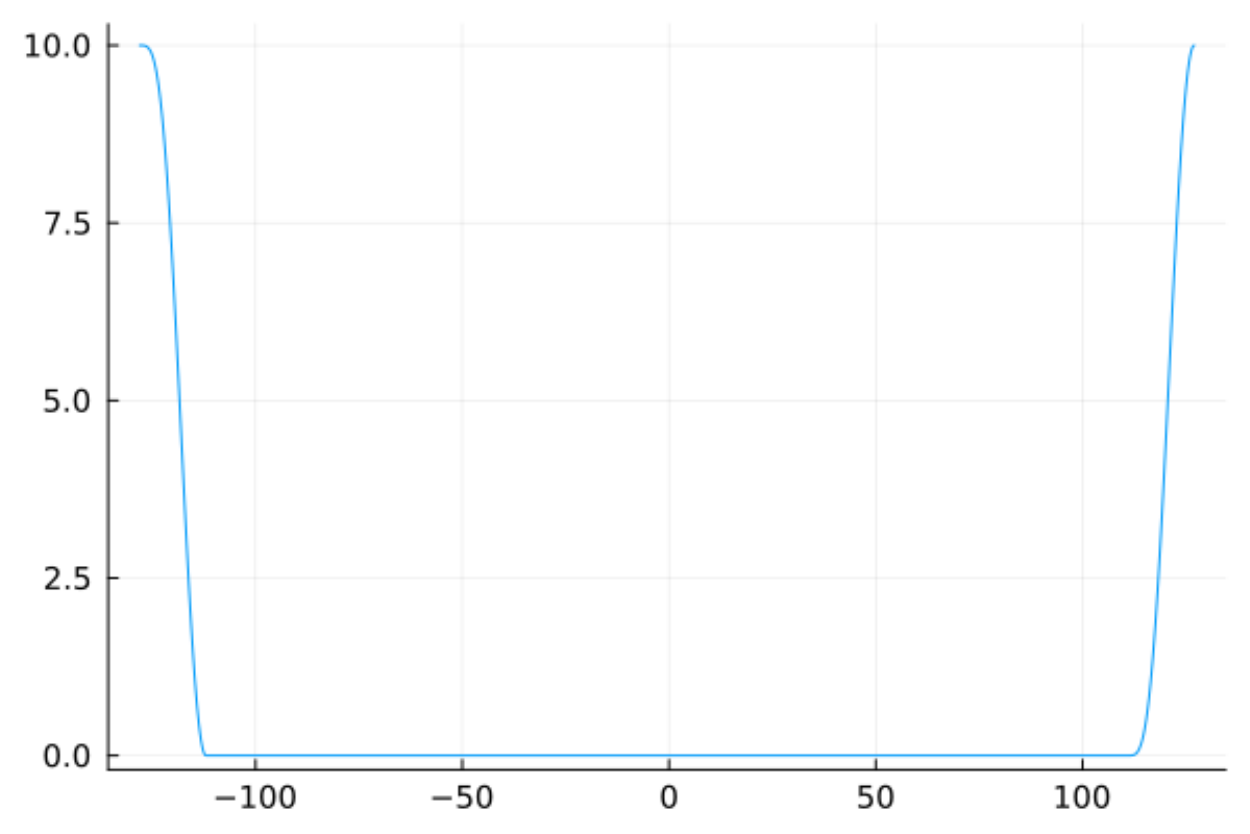}
      \put(5.2,67.5){\tiny $d(x)$}
      \put(100,4){\tiny $x$}
    \end{overpic}
     \begin{figuretext}\label{dampfig}
       The damping function $d(x)$ for $N=128$, $d_0 = 10$, $N_d = \lfloor N/8\rfloor$. 
\end{figuretext}
     \end{center}
\end{figure}

Let $s:\R \to [0,1]$ be continuous step function such that $s(x) = 0$ for $x \leq 0$ and $s(x) = 1$ for $x \geq 1$. Then, for any $d_0 > 0$, 
$$d(x) = d_0 \times \begin{cases}
 1-s\big(\frac{x+N}{N_d}\big) & \text{for $x \in [-N, -N + N_d]$},
	\\
0 & \text{for $x \in (-N + N_d, N-1-N_d)$},
	\\
s\big(\frac{x-N+N_d+1}{N_d}\big) & \text{for $x \in [N -1-N_d, N-1]$},
\end{cases}
$$
is a continuous function of $x \in [-N, N-1]$ such that $d(-N)=d(N-1) = d_0$, $d = 0$ on $(-N + N_d, N-1-N_d)$, and $d(x)$ increases continuously from $0$ to $d_0$ as $x$ increases from $N -1- N_d$ to $N-1$ and as $x$ decreases from $-N + N_d$ to $-N$, see Figure \ref{dampfig}. For the numerical experiments in this paper, we have used $d_0 = 10$, $N_d = \lfloor N/8\rfloor$ where $\lfloor x \rfloor$ is the integer part of $x$, and the following step function which is continuously differentiable on $\R$:
\begin{align}\label{smoothstepdef}
s(x) = \begin{cases} 0 & \text{for $x \leq 0$},
	\\
x^4(x-2)^4 & \text{for $0 \leq x \leq 1$},
	\\
1 & \text{for $x \geq 1$}.
\end{cases}
\end{align}
%Perhaps would be more natural to pick an $s(x)$ such that $s(x)+s(1-x)=1$. For example,
%\begin{align*}
%s(x) = \begin{cases} 0 & \text{for $x \leq 0$},
%	\\
%x^{4}(35-84x+70x^{2}-20x^{3}) & \text{for $0 \leq x \leq 1$},
%	\\
%1 & \text{for $x \geq 1$},
%\end{cases}
%\end{align*}
%is symmetric and is $C^{3}(\R)$. 

We consider the following modification of (\ref{boussinesqsystemFourier}) that includes the additional damping term $-d(j) \hat{U}(j,t)$:
\begin{align}\label{boussinesqsystemFourierdamped}
& \begin{cases}
 \hat{U}_t(j,t) = \frac{\pi i j}{L}\hat{V}(j,t) - d(j) \hat{U}(j,t),
	\\
 \hat{V}_{t}(j,t) = \frac{\pi i j}{L}\hat{U}(j,t) + 2(\widehat{U_{x}} * \hat{U})(j,t) + (\frac{\pi i j}{L})^3 \hat{U}(j,t),
\end{cases} \quad j = -N, \dots, N-1.
\end{align}
One effect of the damping is to smoothen out the Fourier transforms as functions of $j$: the damping suppresses Fourier modes with $j$ close to the endpoints $-N$ and $N-1$, and the closer $j$ gets to these endpoints, the more damped the $j$th mode is.

\subsection{Summary}
To summarize, our numerical scheme for the solution of the initial value problem for (\ref{badboussinesq}) consists of solving the system of $4N$ ODEs (\ref{boussinesqsystemFourierdamped}) with initial conditions (\ref{hatUVinitialconditions}) and $N$ given by (\ref{Ndef}). The system (\ref{boussinesqsystemFourierdamped}) can, for instance, be solved with the help of a Runge--Kutta method. 
The numerical approximation $U(x,t)$ of the solution of the initial value problem for (\ref{badboussinesq}) is then given by (\ref{Udef}).

\section{The small-amplitude long-wave regime}\label{regimesec}
Before applying the numerical scheme of Section \ref{schemesec} to concrete examples, it is useful to elaborate on the regime where (\ref{badboussinesq}) is expected to model water waves. 
The derivation of (\ref{badboussinesq}) relies on two parameters being small: the \textit{amplitude parameter} $\epsilon$ and the \textit{long wavelength}  or \textit{shallowness parameter} $\delta$, see \cite{J1997}. These parameters are defined by $\epsilon = a/h$ and $\delta = h/\lambda$, where $a$ is the typical amplitude, $h$ is the typical depth of the water, and $\lambda$ is the typical wavelength of the surface wave. The limit $\delta \to 0$ is referred to as the long wave or shallow water approximation, while the limit $\epsilon \to 0$ is the small-amplitude approximation. Equation (\ref{badboussinesq}) is a model for long waves of small amplitude, meaning that both $\epsilon$ and $\delta$  are assumed to be small. 

In physical units, (\ref{badboussinesq}) can be written as (see \cite[Eq. (26)]{B1872})
\begin{align}\label{badboussinesqphysical}
\eta_{\tau \tau} - g h \eta_{\xi\xi} - g h \frac{\partial^2}{\partial\xi^2} \bigg(\frac{3}{2 h} \eta^2 + \frac{h^2}{3} \eta_{\xi\xi}\bigg) = 0,
\end{align}
where $g$ is the gravitational acceleration, and $\eta(\xi,\tau)$ is the surface elevation measured with respect to the mean water level at the horizontal position  $\xi$ at time $\tau$. The dimensionless equation (\ref{badboussinesq}) is obtained from (\ref{badboussinesqphysical}) by defining
$$u(x,t) = \frac{3}{2h}\eta(\xi, \tau), \quad x = \xi \frac{\sqrt{3}}{h}, \quad
t = \tau \sqrt{\frac{3g}{h}}.$$ 
Thus, if we denote the typical amplitude and typical wavelength of a solution $u(x,t)$ of (\ref{badboussinesq}) by $A$ and $\Lambda$, respectively, then the parameters $\epsilon$ and $\delta$ are given by $\epsilon = 2A/3$ and $\delta = \sqrt{3}/\Lambda$. If we interpret $\epsilon$ and $\delta$ as being small when they are $< 0.1$, then we should restrict ourselves to solutions of (\ref{badboussinesq}) of amplitude less than $0.15$ and typical wavelength greater than $10 \sqrt{3} \approx 17.3$. The choice we made in Section \ref{schemesec} of setting frequency modes with wave length less than $2\pi$ to zero should be viewed in this light. 

We finally mention that since $g \approx 9.8$ $m/s^2$, the physical time variable $\tau$ is related to $t$  by $\tau \approx 0.184 \sqrt{h}\, t$ seconds. For example, if the mean water depth $h$ is $1$ meter, then $18.4$ seconds have elapsed when $t = 100$.
  
%If $u$ has amplitude $A$, then $\eta$ has amplitude $a = \frac{2h}{3} A$. So $\epsilon = \frac{2h}{3} A/h_0 = \frac{2A}{3}$.
%If $u$ has wavelength $\Lambda$, then $\eta$ has wavelength $\lambda = \frac{h}{\sqrt{3}} \Lambda$. So $\delta = h_0/\lambda = \frac{\sqrt{3}}{\Lambda}$.

%If $t = 1$ and $h = 1$ meter, then $\tau = \sqrt{\frac{h}{3g}} = \sqrt{\frac{1}{3*9.8}} \approx 0.184$ seconds.
%$g = 9.8$ $m/s^2$.
%So if the depth of the water is $1$ meter, then 18.4 seconds have elapsed when $t = 100$.
%If $t = 1$ and $h = 0.5$ meters, then $\tau = \sqrt{\frac{h}{3g}} = \sqrt{\frac{1/2}{3*9.8}} \approx 0.130$ seconds.
%5

%the {\it amplitude parameter} $\epsilon = a/h_0$: the linearised problem; This requires that the amplitude of the surface wave be small; below \cite[(1.62)]{J1997}
%the {\it long wavelength} or {\it shallowness parameter} $\delta$ the long-wave (or shallow-water) problem: the waves are long; that is, of long wavelength (or the water is shallow), in the sense that $\delta = h_0/\lambda$ is small.  below \cite[(1.58)]{J1997}

\section{Comparison with one-soliton solutions}\label{solitonsec}

%\subsection{Comparison with one-solitons}
Equation (\ref{badboussinesq}) admits the following family of exact one-soliton solutions parametrized by $A > 0$ and $x_0 \in \R$:
\begin{align}\label{onesoliton}
u(x,t) = A \sech^2(\sqrt{A/6}(x - x_0 - ct)),
\end{align}
where the wave speed $c$ is given in terms of the amplitude $A$ by  $c = \sqrt{1+ \frac{2A}{3}} > 1$ \cite{H1973} (see also \cite[Appendix A]{CLscatteringsolitons}). 
To test our numerical scheme, we choose $A = 0.05$ (corresponding to a physical wave amplitude that is $3.33$\% of the mean water depth, see Section \ref{regimesec}) and $x_0 = 0$ (so that the peak of the soliton lies initially at $x = 0$).

\begin{figure}
\bigskip\begin{center}
\hspace{-.4cm}
\begin{overpic}[width=.46\textwidth]{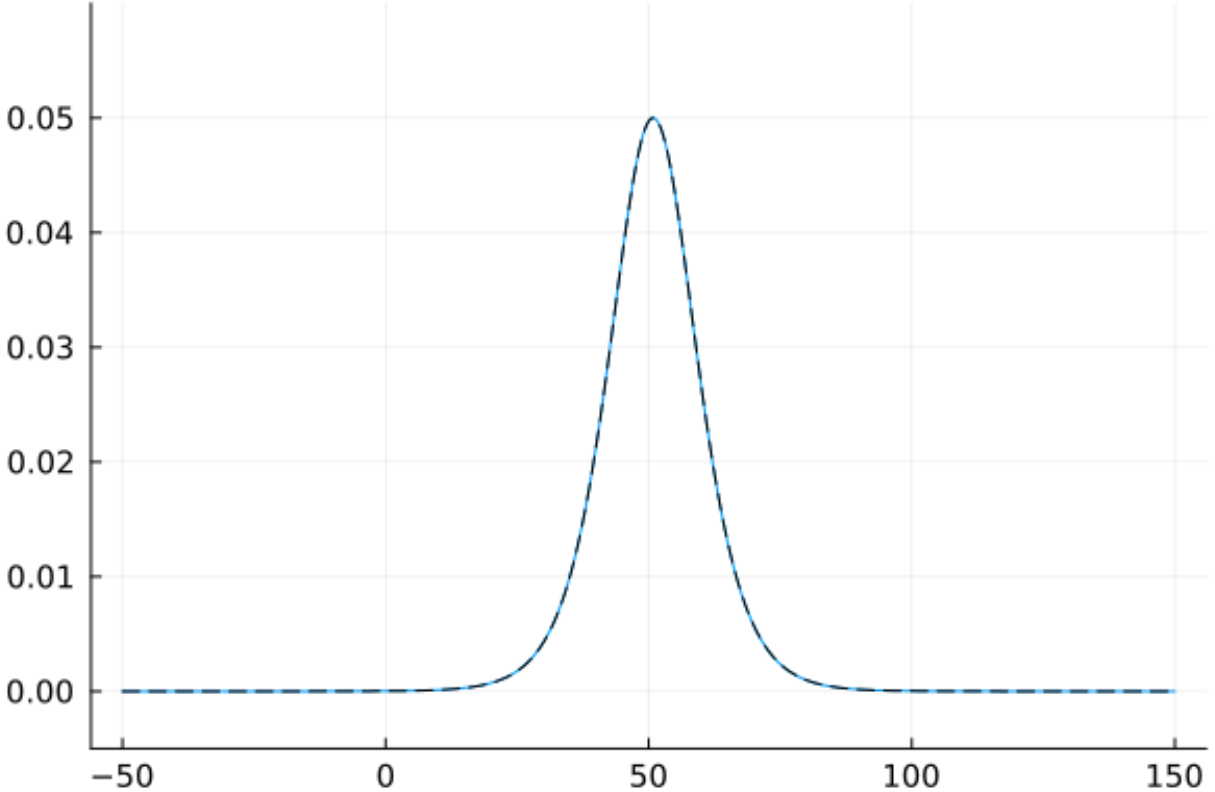}
      \put(-.5,68){\footnotesize $u, U$}
      \put(101.6,2.2){\footnotesize $x$}
    \end{overpic}
    \hspace{0.5cm}
\begin{overpic}[width=.46\textwidth]{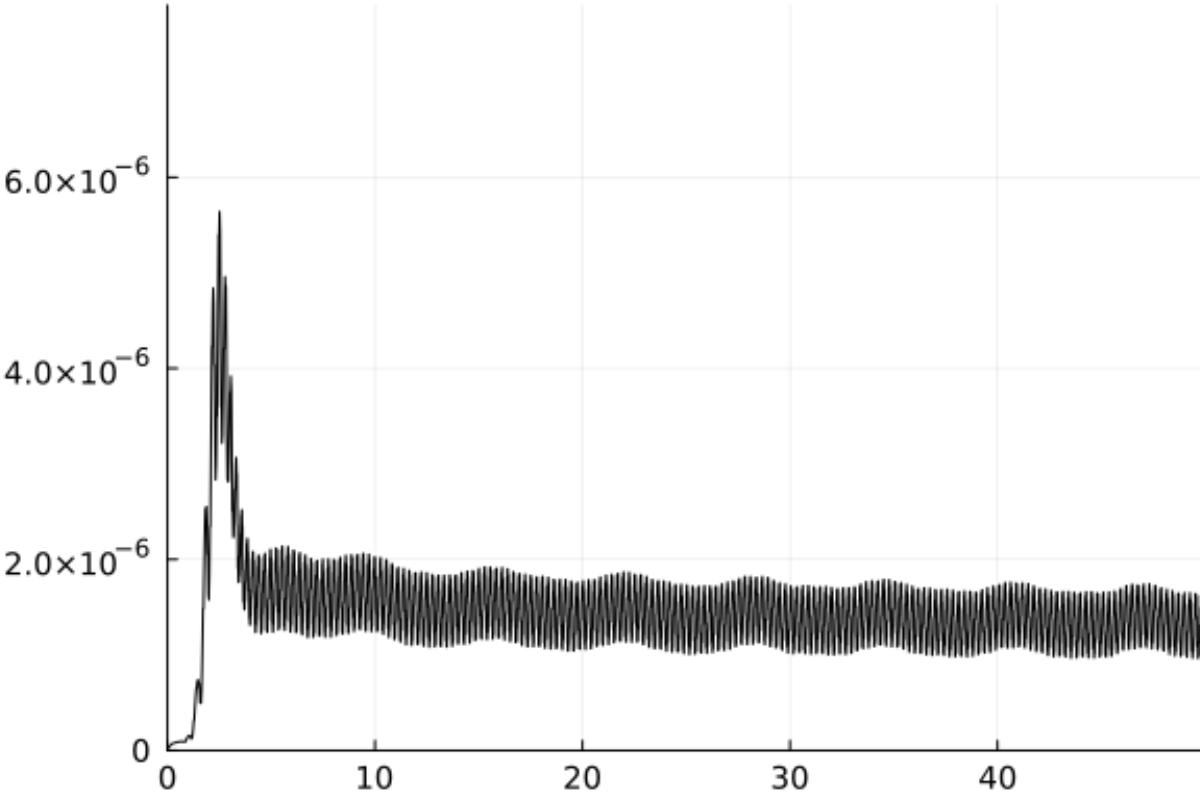}
      \put(11.5,68){\footnotesize $e(t)$}
      \put(101.6,2.2){\footnotesize $t$}
    \end{overpic}
     \begin{figuretext}\label{onesolitonsmallAfig}
       Left: The one-soliton solution $u(x,t)$ (dashed black) for $A = 0.05$ and $x_0 = 0$ together with the numerical approximation $U(x,t)$ (solid blue) obtained with the scheme described in Section \ref{schemesec} at time $t = 50$. (Curves are indistinguishable at this scale.)
       
       Right: The $L^\infty$-error $e(t) = \sup_{x \in [-L,L]}|U(x,t) - u(x,t)|$ as a function of $t \in [0, 50]$.
\end{figuretext}
     \end{center}
\end{figure}

We first consider the soliton for times $t \in [0,50]$. By choosing $L = 200$, we ensure that the wave $u(x,t)$ is negligibly small outside $[-L,L]$ for $t \in [0, 50]$. Applying the numerical scheme described in Section \ref{schemesec} (with $N_d = N/8$, $d_0 = 10$, and $s(x)$ given by (\ref{smoothstepdef})), we obtain an approximation $U(x,t)$ of the one-soliton $u(x,t)$ for $t \in [0, 50]$. Figure \ref{onesolitonsmallAfig} shows $U(x,t)$ and $u(x,t)$ at time $t = 50$, as well as the $L^\infty$-error $e(t)$ defined by\footnote{Numerically, we compute the error $e(t)$ in (\ref{edef}) as $\max_{x \in S} |U(x,t) - u(x,t)|$ where $S$ is a set of $10^5$ equally spaced points ranging from $-L$ to $L$.}
\begin{align}\label{edef}
e(t) := \sup_{x \in [-L,L]}|U(x,t) - u(x,t)|.
\end{align}
As seen in Figure \ref{onesolitonsmallAfig}, the error $e(t)$ reaches its maximum at $t \approx 2.5$ and then decreases with $t$ on the given interval. This is because the numerics leads to an initial drop of the soliton's peak. The peak bounces back up, but $|U(x,t) - u(x,t)|$ remains maximized near $x \approx 0$ for any $t$ on the given interval. For $x$ near the position $ct$ of the soliton's peak, the error $|U(x,t) - u(x,t)|$ grows slightly with $t$, but it is nevertheless smaller than the error near $x = 0$ and therefore not visible in Figure \ref{onesolitonsmallAfig}.

\begin{figure}
\bigskip\begin{center}
\hspace{-.4cm}
\begin{overpic}[width=.46\textwidth]{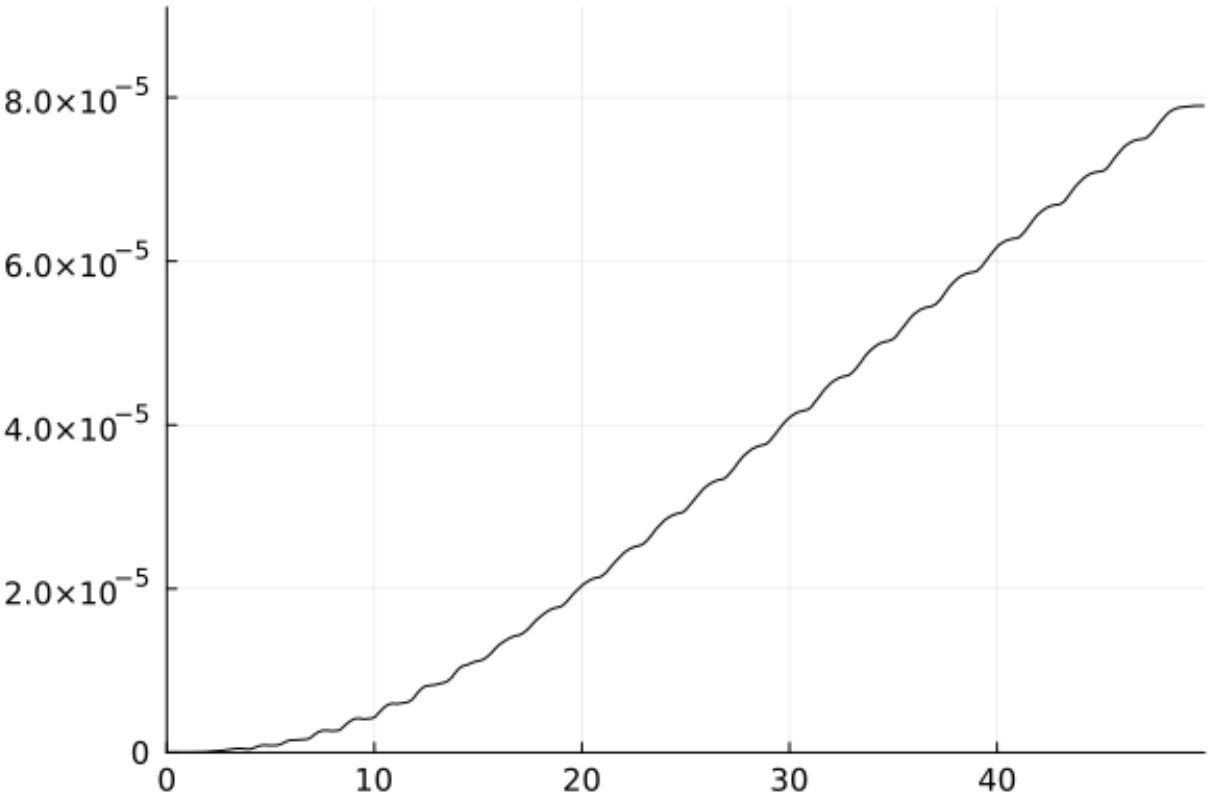}
      \put(11.5,68){\footnotesize $e(t)$}
      \put(101.6,2.2){\footnotesize $t$}
    \end{overpic}
     \begin{figuretext}[No damping gives larger error]\label{onesolitonsmallAnodampingfig}
       The $L^\infty$-error $e(t)$ as a function of $t \in [0, 50]$ for the same simulation as in Figure \ref{onesolitonsmallAfig} except that no damping has been included in the numerical scheme.
\end{figuretext}
     \end{center}
\end{figure}

To demonstrate the effect of the damping introduced in Section \ref{dampingsubsec}, we display in Figure \ref{onesolitonsmallAnodampingfig} the $L^\infty$-error $e(t)$ in the case when no damping has been used in the numerical scheme. In this case, $e(t)$ is more than an order of magnitude larger and grows with $t$.

\begin{figure}
\bigskip\begin{center}
\hspace{-.4cm}
\begin{overpic}[width=.46\textwidth]{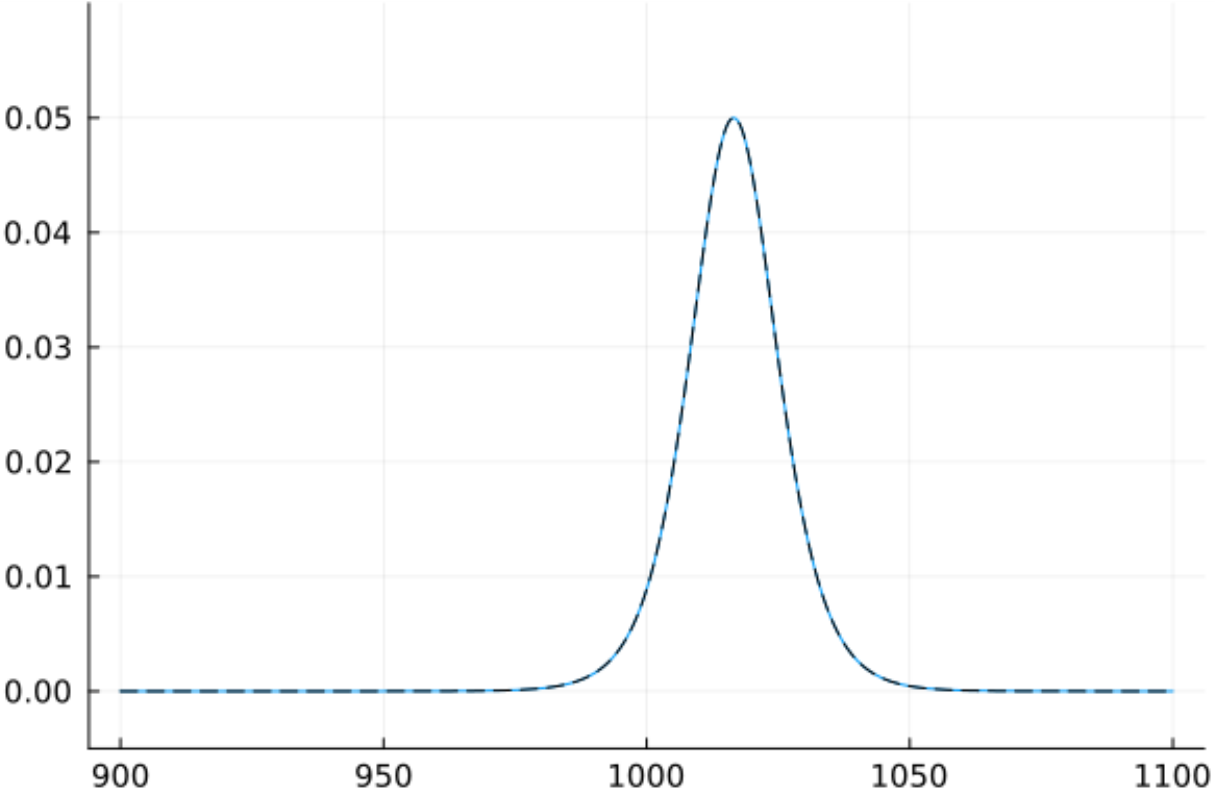}
      \put(-.5,68){\footnotesize $u, U$}
      \put(101.6,2.2){\footnotesize $x$}
    \end{overpic}
    \hspace{0.5cm}
\begin{overpic}[width=.46\textwidth]{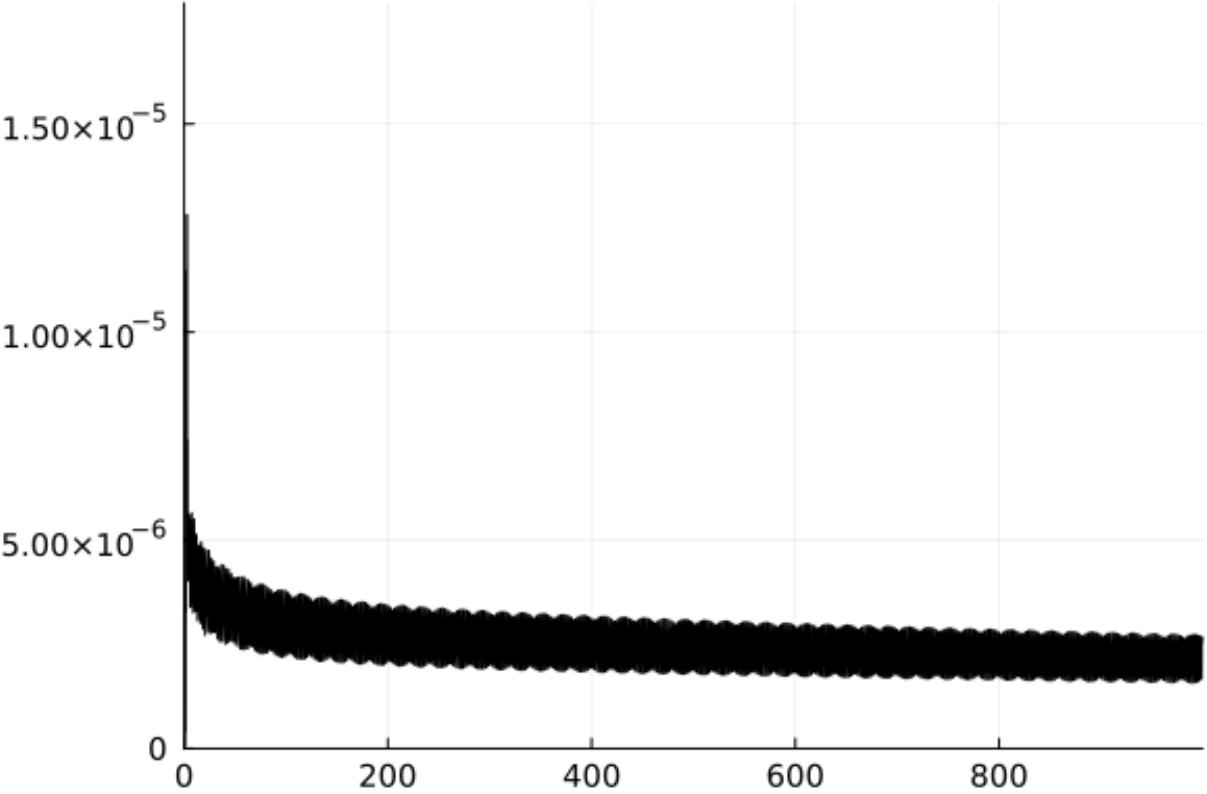}
      \put(12,68){\footnotesize $e(t)$}
      \put(101.6,2.2){\footnotesize $t$}
    \end{overpic}
     \begin{figuretext}\label{onesolitonsmallA1000fig}
       Same as Figure \ref{onesolitonsmallAfig} except that $t = 1000$.
\end{figuretext}
     \end{center}
\end{figure}

To investigate the accuracy of the numerical scheme of Section \ref{schemesec} also on larger time scales, we consider the same soliton solution as above on the time interval $t \in [0, 1000]$. 
Choosing $L = 1200$, we obtain the results shown in Figure \ref{onesolitonsmallA1000fig}. We see that the error remains small for all $t \in [0, 1000]$. The error in Figure \ref{onesolitonsmallA1000fig} is slightly larger than the error in Figure \ref{onesolitonsmallAfig} due to the larger value of $L$.

\subsection{Comparison with other numerical schemes}
\begin{figure}
\bigskip\begin{center}
\hspace{-.4cm}
\begin{overpic}[width=.46\textwidth]{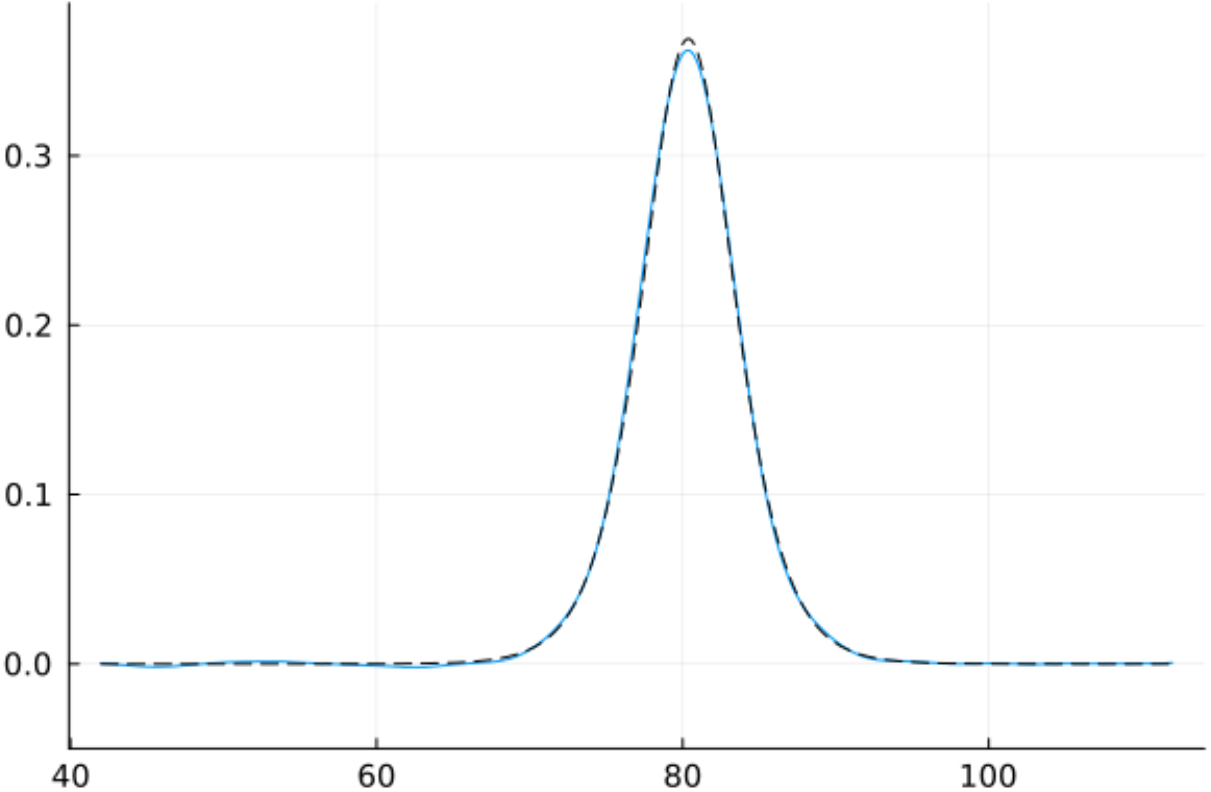}
      \put(-1.5,68){\footnotesize $u, U$}
      \put(101.6,2.6){\footnotesize $x$}
    \end{overpic}
    \hspace{0.5cm}
\begin{overpic}[width=.46\textwidth]{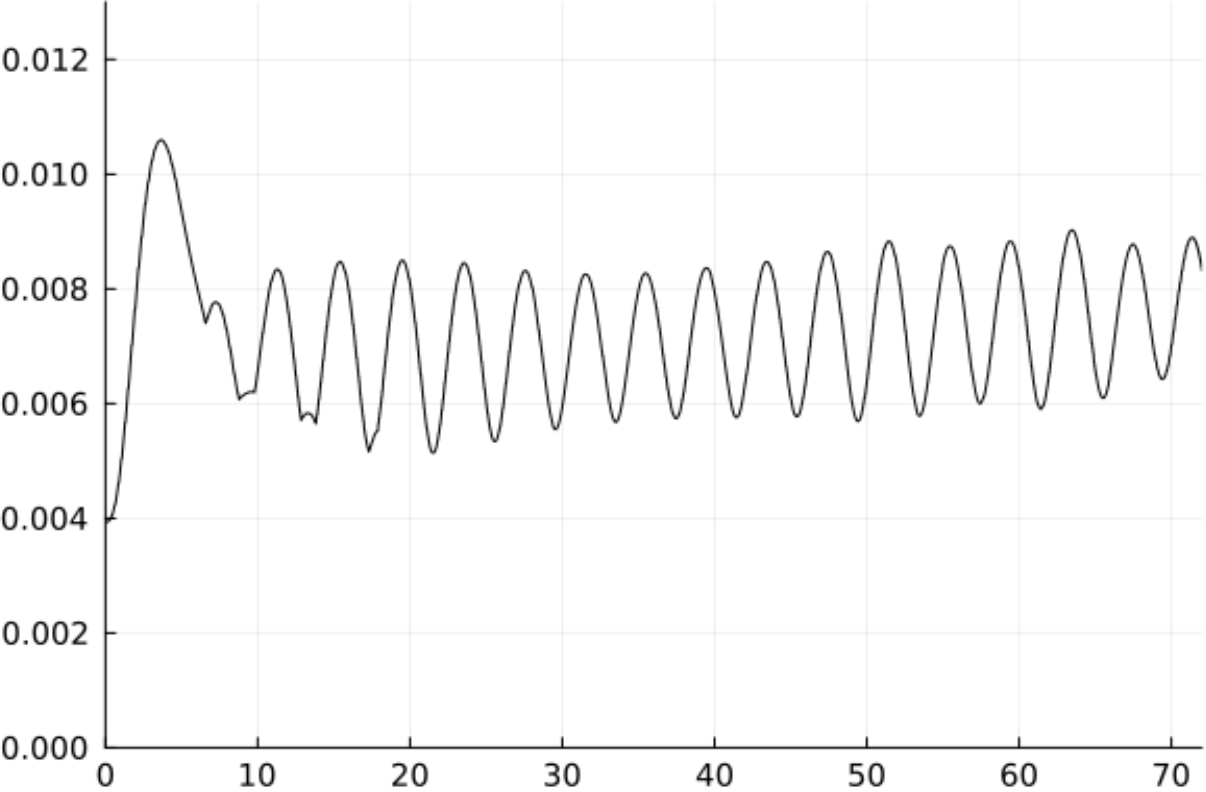}
      \put(5.5,68){\footnotesize $e(t)$}
      \put(101.6,2.2){\footnotesize $t$}
    \end{overpic}
     \begin{figuretext}\label{Bratsosfig}
       Left: The one-soliton solution $u(x,t)$ (dashed black) for $A = 0.369$ and $x_0 = 0$ together with the numerical approximation $U(x,t)$ (solid blue) obtained with the scheme described in Section \ref{schemesec} at time $t = 72$. 
       
       Right: The $L^\infty$-error $e(t) = \sup_{x \in [-L,L]}|U(x,t) - u(x,t)|$ as a function of $t \in [0, 72]$.
\end{figuretext}
     \end{center}
\end{figure}

Other numerical schemes for the solution of (\ref{badboussinesq}) have been proposed and studied in \cite{DH1999, IB2003, BTN2005, B1998, B2007, UEK2021}. The authors of \cite{DH1999} use the one-soliton (\ref{onesoliton}) of amplitude $A = 0.5$ to test their scheme, while in \cite{IB2003, BTN2005, B1998, B2007} the proposed schemes have been applied primarily to the one-soliton of amplitude $A = 0.369$. These amplitudes (corresponding to physical amplitudes that are 33.3\% and 24.6\% of the water's depth, respectively) are outside the regime where we expect (\ref{badboussinesq}) to be an accurate model for water waves, and also outside the range where the numerical scheme of Section \ref{schemesec} is expected to perform well. Nevertheless, we have applied our scheme with $L=200$ to these waves and found good results; see Figure \ref{Bratsosfig} for the results in the case of $A = 0.369$.\footnote{In this case, we have not included any damping in the numerical scheme, because we do not want to suppress more Fourier modes than necessary in these large-amplitude waves.} At $t = 36$ and $t = 72$, the $L^\infty$-error is $0.00839$ and $0.00880$, respectively. The corresponding errors in \cite[Table 1]{B1998} are $0.0358$ and $0.0347$, while the corresponding errors in \cite[Table 9]{UEK2021} are $0.0118$ and $0.0207$. These values and Figure \ref{Bratsosfig} can also be compared with \cite[Figure 2]{B1998}, \cite[Figure 1]{BTN2005}, and \cite[Figure 1b]{B2007}, which show larger errors.

%\cite[Figure 2]{B1998} has error in sup norm about 0.04 for t=72, we get just above 0.01
%About 0.04-0.05 also in \cite[Figure 1]{BTN2005}
%About 0.04 also in \cite[Figure 1b]{B2007}.
%\cite[Table 2]{IB2003} only lists times up until $t = 50$. At $t = 50$, they already have $e \approx 0.074$.
In \cite{B1998, UEK2021}, solitons of other amplitudes are also considered. For the soliton of amplitude $0.1$, \cite[Table 3]{B1998} and \cite[Table 9]{UEK2021} list the $L^\infty$-errors $0.262 \times 10^{-2}$ and $0.147 \times 10^{-3}$ at $t = 36$, respectively, whereas we find $0.212 \times 10^{-4}$. At $t = 72$, \cite[Table 9]{UEK2021} reports the $L^\infty$-error $0.532 \times 10^{-3}$, while our scheme yields $0.204 \times 10^{-4}$. The relatively small errors generated by the scheme of Section \ref{schemesec} are a reflection of the fact that the soliton with $A = 0.1$ falls within the small-amplitude long-wave regime.

In \cite{DH1999, IB2003, BTN2005, B1998, B2007, UEK2021}, the $L^\infty$-errors typically grow relatively quickly as $t$ increases. 
No numerical results have been reported in \cite{DH1999} for $t>8$, in \cite{B1998, BTN2005, B2007} for $t>72$ and in \cite{IB2003, UEK2021} for $t>120$. 
A strength of the scheme presented in this paper is that it performs very well for moderate and large values of $t$. 
In these regimes, our scheme yields significantly smaller $L^\infty$-errors than those reported in \cite{IB2003, BTN2005, B1998, B2007, UEK2021}. Our scheme performs less well for small values of $t$, where the loss of high-frequency modes leads to a temporarily increased error (as in Figure \ref{onesolitonsmallAfig}).

\section{Comparison with asymptotic formulas}\label{asymptoticsec}
Asymptotic formulas for the long-time behavior of the solution of (\ref{badboussinesq}) were presented in \cite{CLmain}. Ten main asymptotic sectors were identified and in each sector an exact expression for the leading asymptotic term together with a precise error estimate was provided. The formulas are valid whenever the initial data has sufficient smoothness and decay and is such that its ``nonlinear high-frequency modes'' vanish (more precisely, a certain spectral function, $r_1(k)$, must vanish on the vertical segment from $0$ to $i$ in the spectral complex plane, see \cite[Assumption 2.11]{CLmain}). As discussed in the introduction, this implies that the associated solution of (\ref{badboussinesq}) exists for all times \cite[Theorem 2.12]{CLmain}. 

In this section, we apply the numerical scheme of Section \ref{schemesec} to various choices of the initial data and compare the approximate solutions it generates with the asymptotic formulas of \cite{CLmain}. Since the asymptotic formulas we use are only valid up to an error of order $O(t^{-1} \ln t)$ 
or $O(t^{-1/2})$ depending on the sector, 
we cannot expect complete agreement. However, as $t$ increases we expect the numerical approximation to approach the leading asymptotic term, and this is indeed what is found.

\subsection{Gaussian initial data}\label{onehumpsubsec}
\begin{figure}
\bigskip\begin{center}
\hspace{-.4cm}
\begin{overpic}[width=.46\textwidth]{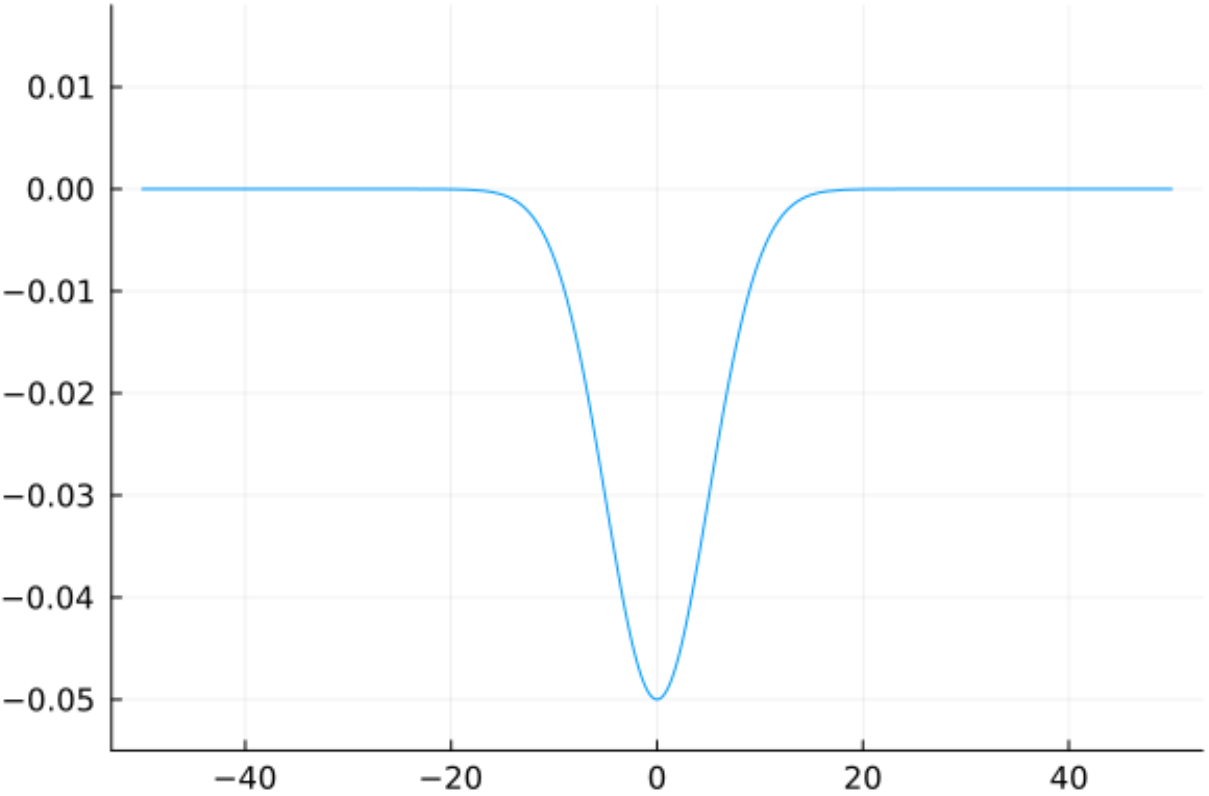}
      \put(8,68){\footnotesize $U$}
      \put(101.6,2.6){\footnotesize $x$}
       \put(49,54){\footnotesize $t = 0$}
    \end{overpic}
    \hspace{0.5cm}
\begin{overpic}[width=.46\textwidth]{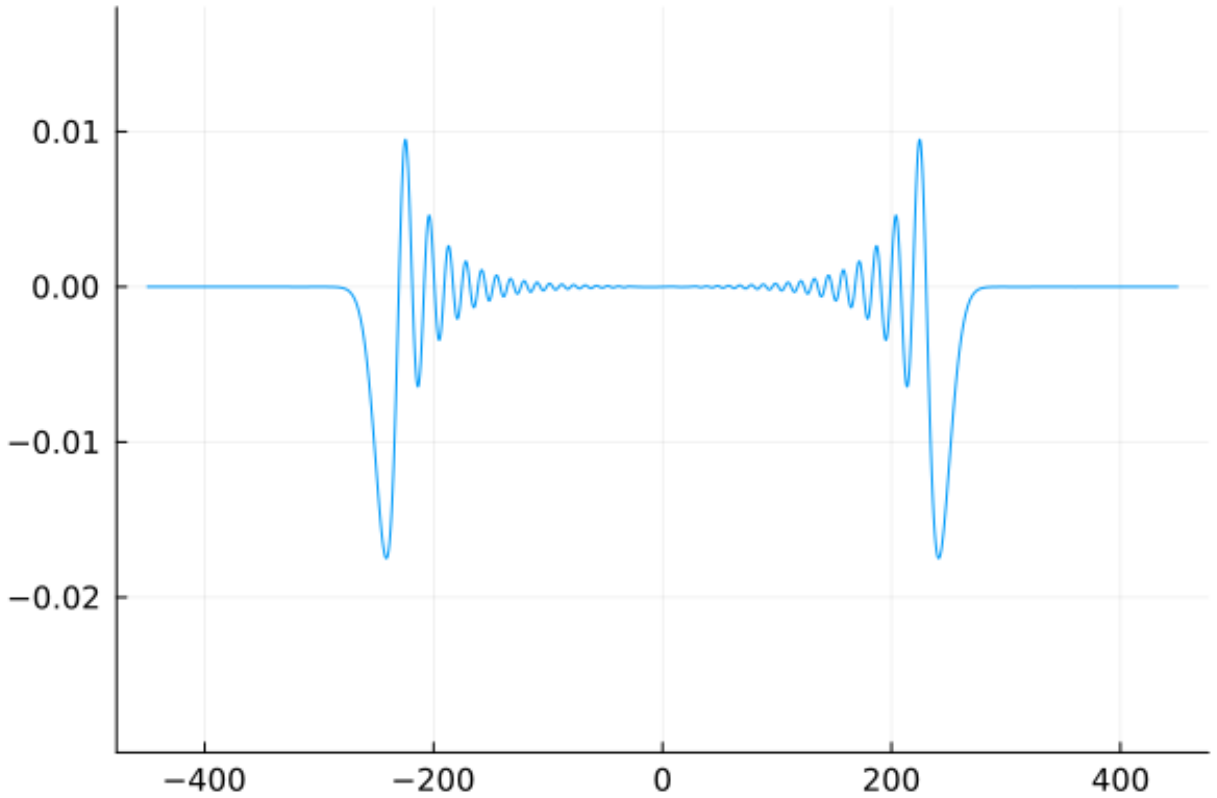}
      \put(8,68){\footnotesize $U$}
      \put(101.6,2.8){\footnotesize $x$}
       \put(47,48){\footnotesize $t = 250$}
    \end{overpic}
    \\
    \vspace{0.5cm}
    \hspace{-0.38cm}
\begin{overpic}[width=.46\textwidth]{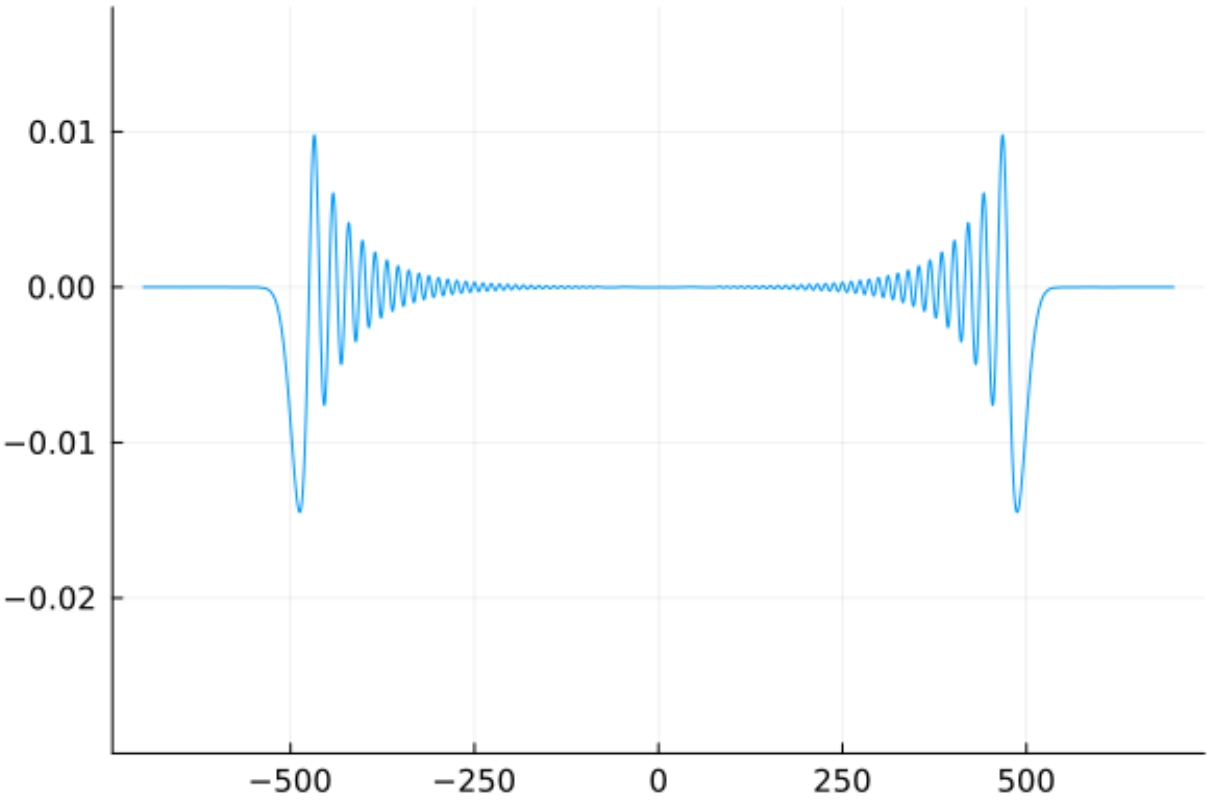}
      \put(8,68){\footnotesize $U$}
      \put(101.6,2.8){\footnotesize $x$}
       \put(47,48){\footnotesize $t = 500$}
    \end{overpic}
        \hspace{0.54cm}
\begin{overpic}[width=.46\textwidth]{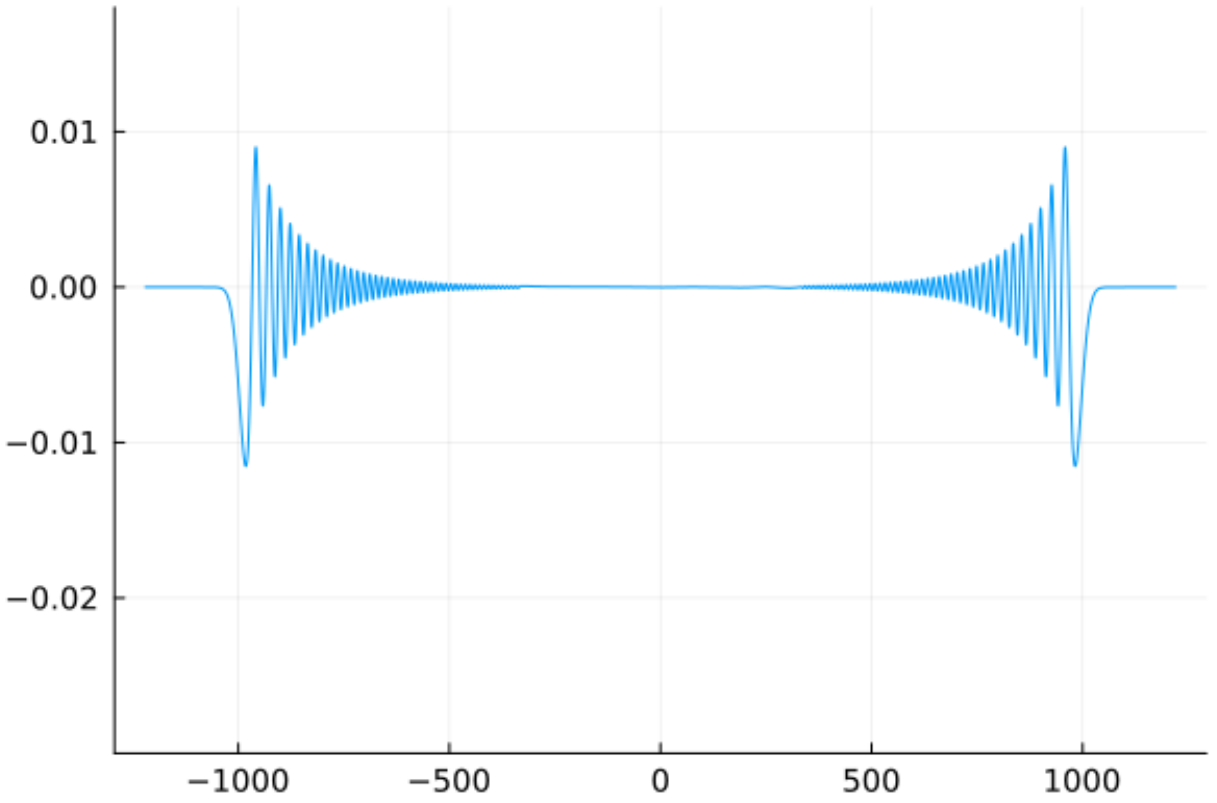}
      \put(8,68){\footnotesize $U$}
      \put(101.6,2.8){\footnotesize $x$}
       \put(47,48){\footnotesize $t = 1000$}
    \end{overpic}
     \begin{figuretext}\label{onehumpfig1}
       The approximate solution $U(x,t)$ of (\ref{badboussinesq}) obtained by the scheme of Section \ref{schemesec} for the Gaussian initial data (\ref{Gaussian}) at times $t = 0$, $t=250$, $t=500$, and $t=1000$. 
\end{figuretext}
     \end{center}
\end{figure}

We consider the Gaussian initial data 
\begin{align}\label{Gaussian}
u(x,0) = -0.05 e^{-0.02 x^2}, \qquad u_t(x,0) = 0.
\end{align}
Applying the scheme of Section \ref{schemesec} with this initial data, we obtain the approximate solution $U(x,t)$ shown in Figure \ref{onehumpfig1}. 

\subsubsection{Comparison with the leading asymptotic term}\label{comparisonsubsubsec}
We focus on the asymptotic sectors denoted by IV and V in \cite{CLmain}, which are given by $\frac{1}{\sqrt{3}} < \zeta < 1$ and $0 < \zeta < \frac{1}{\sqrt{3}}$, respectively, where $\zeta := x/t$. According to \cite{CLmain}, the solution $u(x,t)$ of (\ref{badboussinesq}) with initial data (\ref{Gaussian}) obeys in Sectors IV and V the asymptotic estimate\footnote{Strictly speaking, the asymptotic formula (\ref{uasymptotics}) only holds for initial data whose nonlinear high-frequency Fourier modes vanish (in a sense made precise in \cite{CLmain}). Since the initial data (\ref{Gaussian}) is very slowly varying, the nonlinear high-frequency Fourier modes are almost, but not exactly, zero. The estimate (\ref{uasymptotics}) is therefore only true up to a small error. We will ignore this error in what follows.}
\begin{align}\label{uasymptotics}
u(x,t) = u_{\lead}(x,t) + O\bigg(\frac{\ln t}{t}\bigg) \qquad \text{as $t \to \infty$},
\end{align}
where the leading asymptotic term $u_{\lead}(x,t)$ has the form
\begin{align}\label{uleaddef}
u_{\lead}(x,t) = \begin{cases} 
\frac{A_{1}(\zeta)}{\sqrt{t}}   \cos \alpha_{1}(\zeta,t) + \frac{A_{2}(\zeta)}{\sqrt{t}} \cos \alpha_{2}(\zeta,t) &\text{in Sector IV (i.e., for $\frac{1}{\sqrt{3}} < \zeta < 1)$}, 
	\\
\frac{\tilde{A}_{1}(\zeta)}{\sqrt{t}}   \cos \tilde{\alpha}_{1}(\zeta,t) +\frac{A_{2}(\zeta)}{\sqrt{t}} \cos \tilde{\alpha}_{2}(\zeta,t) &\text{in Sector V (i.e., for $0 < \zeta < \frac{1}{\sqrt{3}}$)}.
\end{cases}
\end{align}
Note that the amplitude function $A_2$ is the same in both sectors. By \cite[Theorem 2.14]{CLmain}, the formula (\ref{uasymptotics}) is valid uniformly for $\zeta$ in compact subsets of $(0, \frac{1}{\sqrt{3}}) \cup (\frac{1}{\sqrt{3}},1)$. 

The functions $A_1, A_2, \alpha_1, \alpha_2, \tilde{A}_1, \tilde{\alpha}_1, \tilde{\alpha}_2$ in (\ref{uleaddef}) are defined in terms of the initial data $u_0(x) = u(x,0)$ and $u_1(x) = u_t(x,0)$ as follows. Let $\omega := e^{\frac{2\pi i}{3}}$, $v_0(x) := \int_{-\infty}^x u_1(x')dx'$, and
\begin{align}\label{lmexpressions intro}
& l_{j}(k) := i \frac{\omega^{j}k + (\omega^{j}k)^{-1}}{2\sqrt{3}}, \qquad j = 1,2,3, \quad k \in \C\setminus \{0\}.
\end{align}
Introduce $P(k)$ and $\mathsf{U}(x,k)$ by
\begin{align*}
P(k) = \begin{pmatrix}
1 & 1 & 1  \\
l_{1}(k) & l_{2}(k) & l_{3}(k) \\
l_{1}(k)^{2} & l_{2}(k)^{2} & l_{3}(k)^{2}
\end{pmatrix}, \quad \mathsf{U}(x,k) = P(k)^{-1} \begin{pmatrix}
0 & 0 & 0 \\
0 & 0 & 0 \\
-\frac{u_{0x}}{4}-\frac{iv_{0}}{4\sqrt{3}} & -\frac{u_{0}}{2} & 0
\end{pmatrix} P(k),
\end{align*} 
and let $X(x,k)$ be the unique solution of the linear integral equation
\begin{align*}  
 & X(x,k) = I - \int_x^{\infty} e^{(x-x')\mathcal{L}(k)} (\mathsf{U}X)(x',k) e^{-(x-x')\mathcal{L}(k)} dx',
\end{align*}
where $\mathcal{L} := \diag(l_1 , l_2 , l_3)$. Define 
\begin{align}\label{sdef}
& s(k) := I - \int_\R e^{-x \mathcal{L}(k)}(\mathsf{U}X)(x,k)e^{x \mathcal{L}(k)}dx,
\qquad r_1(k) := \frac{s_{12}(k)}{s_{11}(k)},
\end{align}
where $s_{ij}(k)$ is the $(ij)$th entry of the $3 \times 3$-matrix $s(k)$.
The amplitude $A_2(\zeta)$ is given in terms of
$$\hat{\nu}_2(k_{2}) := - \frac{1}{2\pi} \ln\bigg(\frac{(1+\tilde{r}(\omega^2 k_{2}) |r_1(\omega^2 k_{2})|^2) |s_{11}(\omega^2 k_{2})|^2}{|s_{11}(\omega k_{2})|^2}\bigg) \geq 0$$
by
$$A_2(\zeta) := \frac{-4\sqrt{3}\sqrt{\hat{\nu}_2(k_2)}|\tilde{r}(\frac{1}{k_{2}})|^{\frac{1}{2}}\im k_2}{-i\omega^{2} k_{2} z_{2,\star}}\sin(\arg(\omega^{2} k_{2})), $$
where $\tilde{r}(k):=\frac{\omega^{2}-k^{2}}{1-\omega^{2}k^{2}}$ and
\begin{align} \label{k2def}
& k_{2} = k_{2}(\zeta) := \frac{1}{4}\bigg( \zeta - \sqrt{8+\zeta^{2}} - i \sqrt{2}\sqrt{4-\zeta^{2}+\zeta\sqrt{8+\zeta^{2}}} \bigg) 
	\\ \nonumber
& z_{2,\star} = z_{2,\star}(\zeta) := \sqrt{2}e^{\frac{\pi i}{4}} \sqrt{-\omega^{2} \frac{4-3k_{2} \zeta - k_{2}^{3} \zeta}{4k_{2}^{4}}}, \qquad -i\omega^{2} k_{2}z_{2,\star}>0.
\end{align}
We refer to \cite{CLmain} for expressions for the functions $A_1, \alpha_1, \alpha_2, \tilde{A}_1, \tilde{\alpha}_1, \tilde{\alpha}_2$. The expressions for the amplitudes $A_1$ and $\tilde{A}_1$ are similar to the expression for $A_2$ given above; the expressions for $\alpha_1, \alpha_2, \tilde{\alpha}_1, \tilde{\alpha}_2$ are more complicated, but just like $A_1, \tilde{A}_1, A_2$ they are expressed in terms of the initial data via linear operations that are numerically accessible. 
We have used the Julia package ``ApproxFun" developed by S. Olver\footnote{See \texttt{https://github.com/JuliaApproximation/ApproxFun.jl} and \cite{OT2013}.} to compute $u_{\lead}(x,t)$ numerically.

\begin{figure}
\bigskip\begin{center}
\hspace{-.4cm}
\begin{overpic}[width=.46\textwidth]{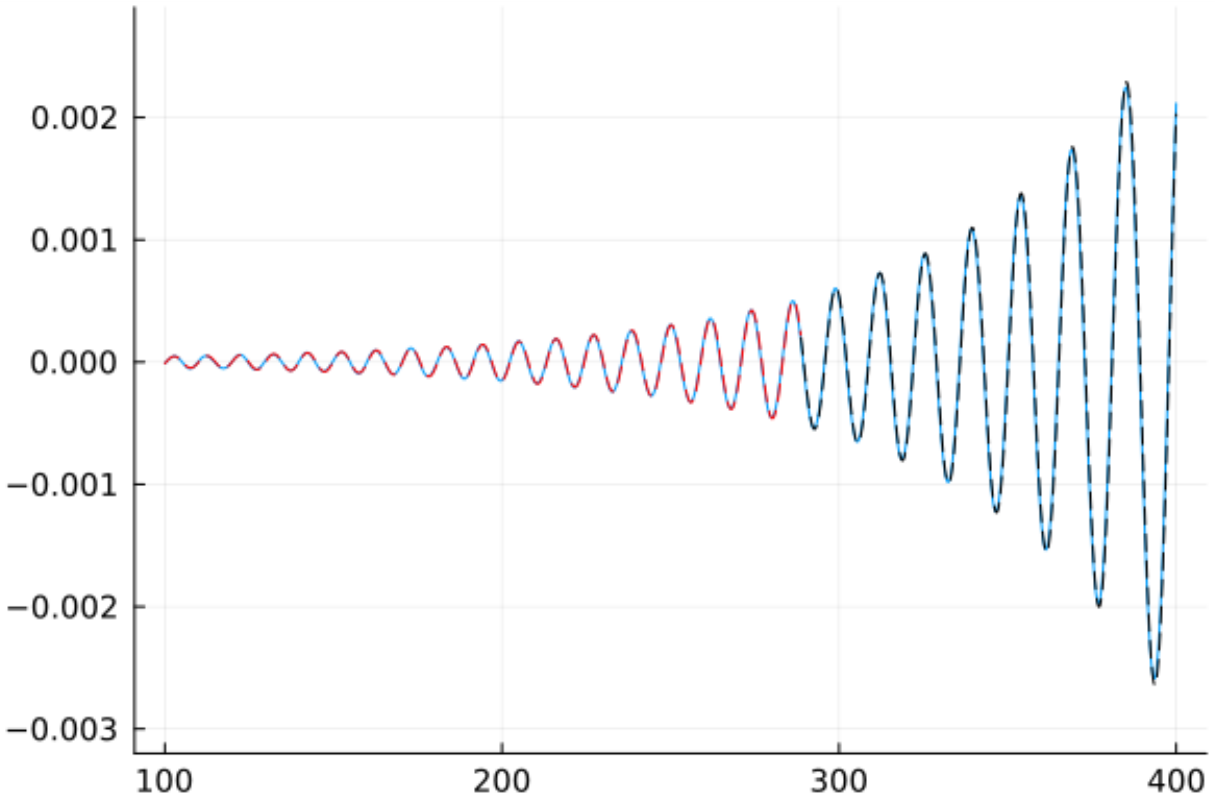}
      \put(4,68){\footnotesize $u_{\lead}, U$}
      \put(101.6,2.2){\footnotesize $x$}
        \put(46,48){\footnotesize $t = 500$}
   \end{overpic}
    \hspace{0.5cm}
\begin{overpic}[width=.46\textwidth]{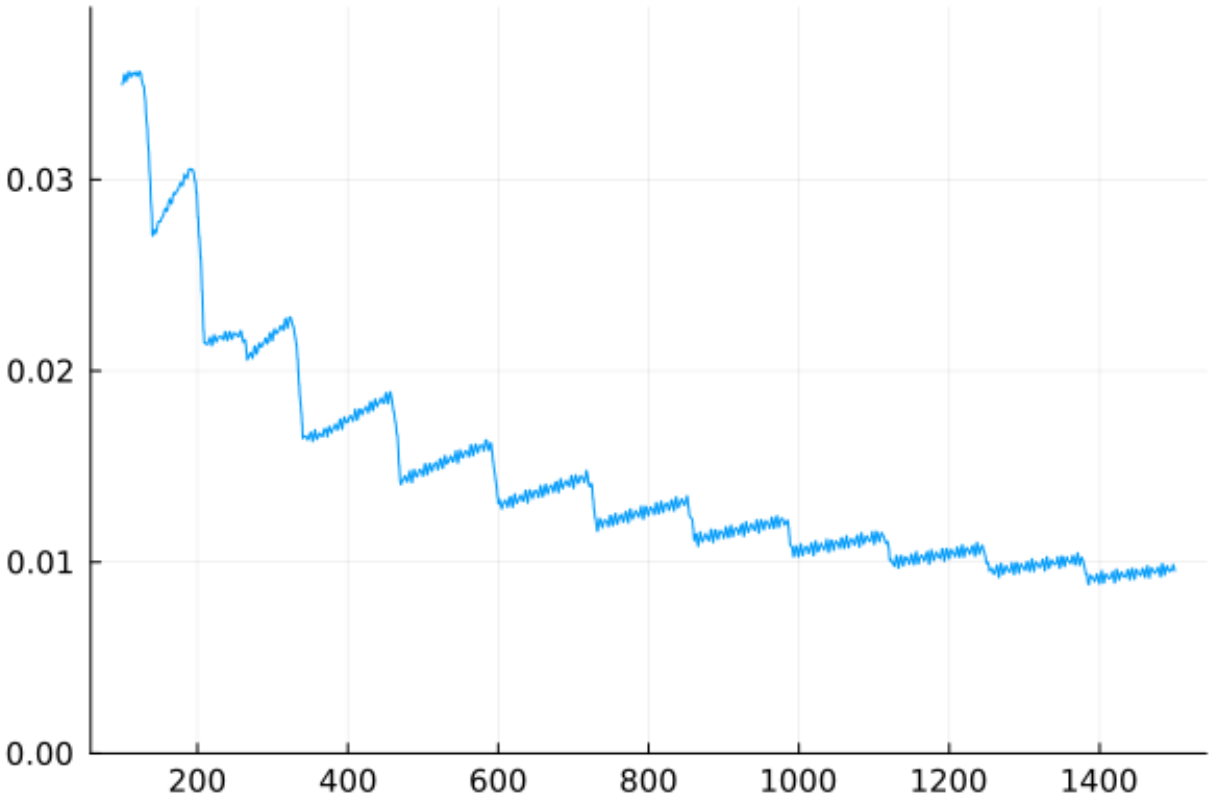}
      \put(0,68){\footnotesize $\frac{t}{\ln t}e(t)$}
      \put(101.6,2.2){\footnotesize $t$}
    \end{overpic}
     \begin{figuretext}\label{onehumpfig2}
       Left: The leading asymptotic term $u_{\lead}(x,t)$ (dashed black in Sector IV and dashed red in Sector V) together with the numerical approximation $U(x,t)$ (solid blue) for the Gaussian initial data (\ref{Gaussian}) obtained with the scheme of Section \ref{schemesec} at time $t = 500$.
       
       Right: The graph of $t \mapsto \frac{t}{\ln t} e(t)$ for $t\in [0,1500]$, where $e(t) = \sup_{x \in [0.2t, 0.8t]} |U(x,t) - u_{\lead}(x,t)|$.
\end{figuretext}
     \end{center}
\end{figure}

In what follows, we compare $u_{\lead}(x,t)$ with the numerical approximation $U(x,t)$ obtained by the scheme of Section \ref{schemesec} for the initial data (\ref{Gaussian}). The left half of Figure \ref{onehumpfig2} shows $u_{\lead}(x,t)$ and $U(x,t)$ for $x \in [0.2t, 0.8t]$ at $t = 500$. To test the accuracy of the numerical scheme, we consider the difference $U(x,t) - u_{\lead}(x,t)$ in the asymptotic sector $\zeta \in [0.2, 0.8]$. From the asymptotic formula (\ref{uasymptotics}), we know that, for any compact subset $\mathcal{I}$ of $(0, \frac{1}{\sqrt{3}}) \cup (\frac{1}{\sqrt{3}},1)$, $\frac{t}{\ln t}\sup_{x/t \in \mathcal{I}}|u(x,t) - u_{\lead}(x,t)|$ remains bounded as $t \to \infty$. In fact, we expect the formula (\ref{uasymptotics}) to be uniform also near $1/\sqrt{3}$, see Section \ref{conjecturesubsec}. In the right half of Figure \ref{onehumpfig2}, we therefore plot $\frac{t}{\ln t} e(t)$ where $e(t) = \sup_{x \in [0.2t, 0.8t]} |U(x,t) - u_{\lead}(x,t)|$ is the $L^\infty$-error on the interval $[0.2t, 0.8t]$. The function $\frac{t}{\ln t} e(t)$ shows no tendency to grow, which is evidence for the accuracy of the numerical scheme of Section \ref{schemesec}.

\begin{figure}
\bigskip\begin{center}
\hspace{-.4cm}
\begin{overpic}[width=.46\textwidth]{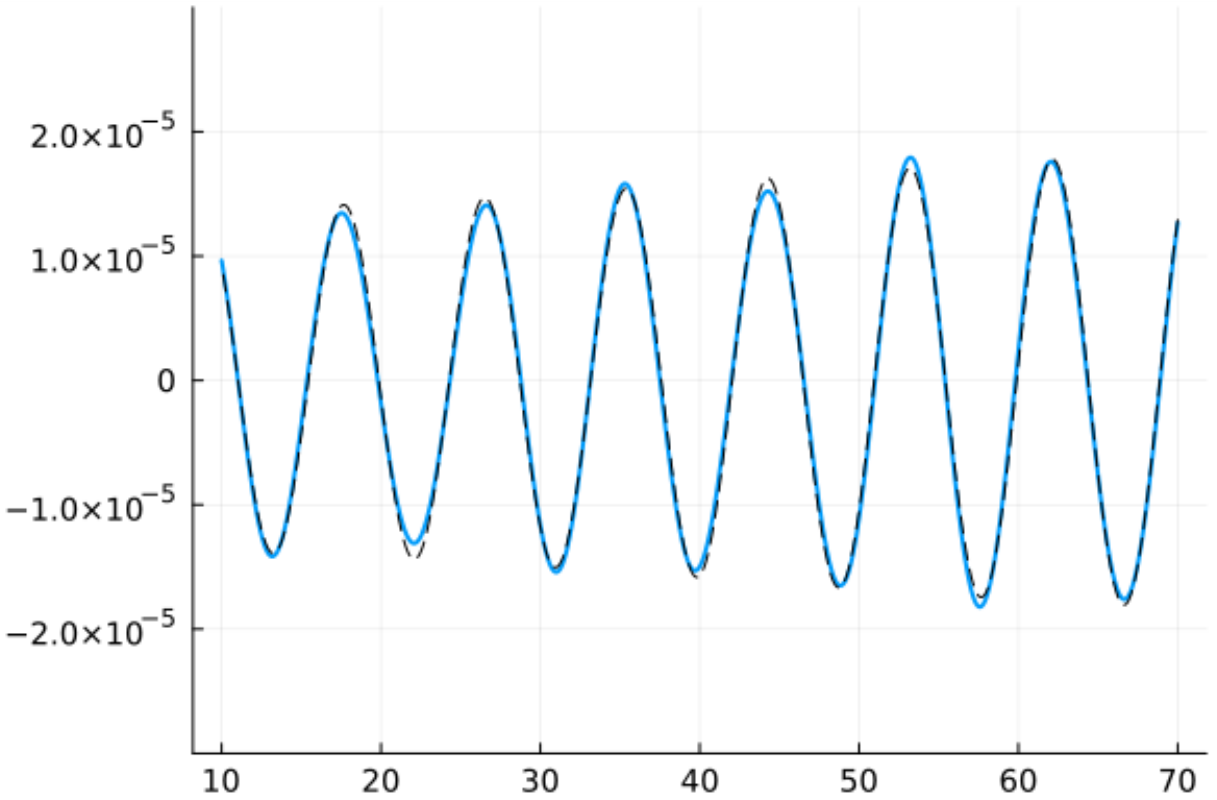}
      \put(5,68){\footnotesize $u_{\lead}, U$}
      \put(101.6,2.2){\footnotesize $x$}
        \put(50,58){\footnotesize $t = 1000$}
    \end{overpic}
    \hspace{0.5cm}
\begin{overpic}[width=.46\textwidth]{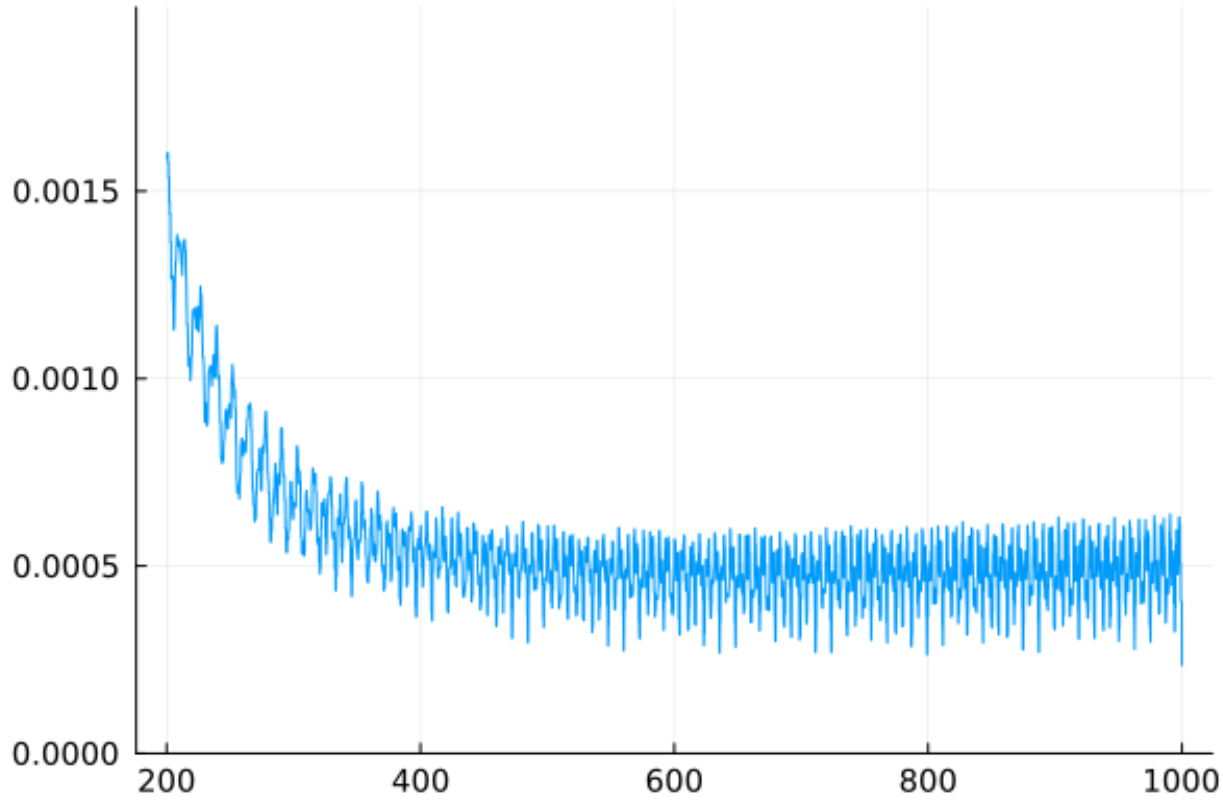}
      \put(4,68){\footnotesize $\frac{t}{\ln t}e(t)$}
      \put(101.6,2.2){\footnotesize $t$}
    \end{overpic}
     \begin{figuretext}\label{interactionfig}
       Left: The leading asymptotic term $u_{\lead}(x,t)$ (dashed black) together with the numerical approximation $U(x,t)$ (solid blue) for the Gaussian initial data (\ref{Gaussian}) obtained with the scheme of Section \ref{schemesec} at time $t = 1000$. 
       
       Right: The graph of $t \mapsto \frac{t}{\ln t} e(t)$ where $e(t) = \sup_{x \in [10,70]}|U(x,t) - u_{\lead}(x,t)|$.
\end{figuretext}
     \end{center}
\end{figure}

\begin{figure}
\bigskip\begin{center}
\hspace{-.4cm}
\begin{overpic}[width=.46\textwidth]{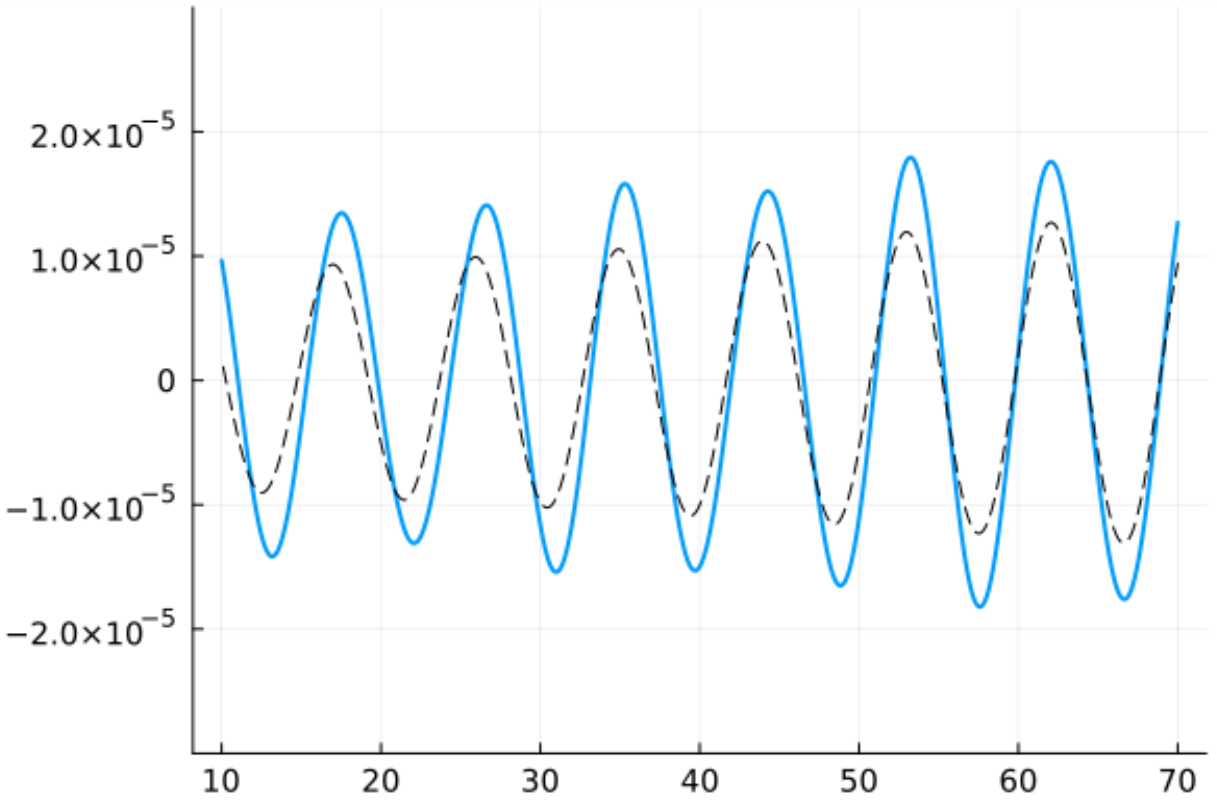}
      \put(6,68){\footnotesize $u_{\lead, r}, U$}
      \put(101.6,2.2){\footnotesize $x$}
        \put(50,58){\footnotesize $t = 1000$}
    \end{overpic}
    \hspace{0.5cm}
\begin{overpic}[width=.46\textwidth]{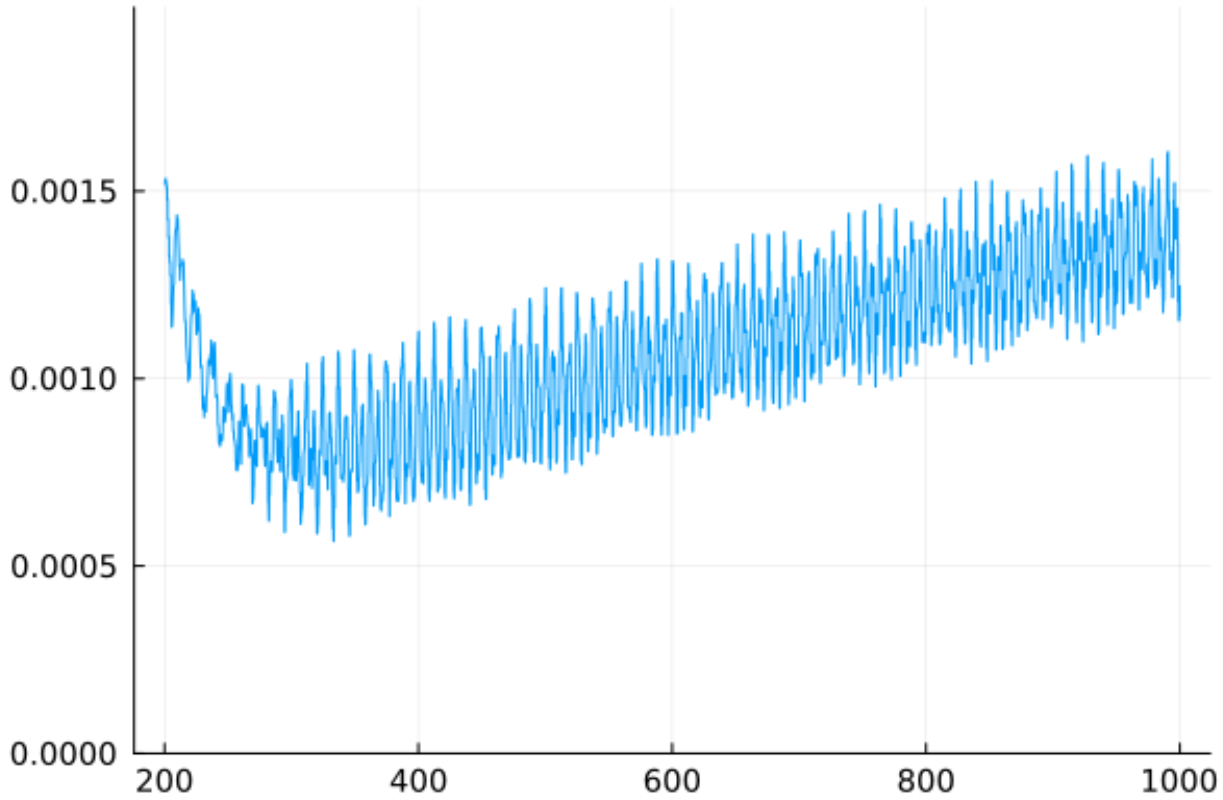}
      \put(3,68){\footnotesize $\frac{t}{\ln t}e_r(t)$}
      \put(101.6,2.2){\footnotesize $t$}
    \end{overpic}
     \begin{figuretext}[No left-moving wave gives large error] \label{interactiononlyrightmovingfig}
       Left: The right-moving wave $u_{\lead,r}(x,t)$ (dashed black) together with the numerical approximation $U(x,t)$ (solid blue) for the initial data (\ref{Gaussian}) obtained with the scheme of Section \ref{schemesec} at time $t = 1000$. 
       
       Right: The graph of $t \mapsto \frac{t}{\ln t} e_r(t)$ where $e_r(t) = \sup_{x \in [10,70]}|U(x,t) - u_{\lead,r}(x,t)|$.
\end{figuretext}
     \end{center}
\end{figure}

\subsubsection{Interaction of right- and left-moving waves}
One of the most interesting aspects of the ``bad'' Boussinesq equation, which sets it apart from for example the KdV equation, is that it is bidirectional. This means that it models waves propagating in both the positive and negative $x$-directions, see e.g. \cite[Section 3.2.5]{J1997}. 
In the asymptotic formula (\ref{uasymptotics}), the bidirectionality manifests itself in that $u_{\lead} = u_{\lead,r} + u_{\lead,l}$ is the sum of a right-moving wave $u_{\lead,r}$ and a left-moving wave $u_{\lead,l}$:
$$u_{\lead}(x,t) = u_{\lead,r}(x,t) + u_{\lead,l}(x,t)$$
where
\begin{align*}
& u_{\lead,r}(x,t) = \frac{A_{1}(\zeta)}{\sqrt{t}}   \cos \alpha_{1}(\zeta,t) \quad \text{and} \quad 
u_{\lead,l}(x,t) =  \frac{A_{2}(\zeta)}{\sqrt{t}} \cos \alpha_{2}(\zeta,t) \quad \text{in Sector IV}; 
	\\
& u_{\lead,r}(x,t) = \frac{\tilde{A}_{1}(\zeta)}{\sqrt{t}}   \cos \tilde{\alpha}_{1}(\zeta,t)\quad \text{and} \quad u_{\lead,l}(x,t) =  \frac{A_{2}(\zeta)}{\sqrt{t}} \cos \tilde{\alpha}_{2}(\zeta,t) \quad \text{in Sector V}.
\end{align*}
Since the bidirectionality is a prominent feature of (\ref{badboussinesq}), it is interesting to determine whether the numerical scheme of Section \ref{schemesec} is accurate enough to resolve the right- and left-moving components of $u_{\lead}$.

Figure \ref{interactionfig} shows $u_{\lead}(x,t)$ and the numerical approximation $U(x,t)$ for $x \in [10,70]$ at time $t = 1000$ for the Gaussian initial data (\ref{Gaussian}). Figure \ref{interactionfig} also shows $\frac{t}{\ln t} e(t)$ where $e(t) = \sup_{x \in [10,70]}|U(x,t) - u_{\lead}(x,t)|$ is the $L^\infty$-error on $[10,70]$. Figure \ref{interactionfig} suggests that, up to small oscillations,  $\frac{t}{\ln t} e(t)$ approaches a constant for large $t$, in agreement with (\ref{uasymptotics}).

%The left-moving wave $u_{\lead,l}$ is generally quite small for  large $x \geq 0$. To observe it, we consider the solution for $x \in [10,70]$. 

To observe the left-moving wave $u_{\lead,l}(x,t)$, we plot in Figure \ref{interactiononlyrightmovingfig} the same graphs as in Figure \ref{interactionfig} except that $u_{\lead}$ has now been replaced with $u_{\lead,r}$, i.e., only the right-moving wave $u_{\lead,r}$ has been included in the asymptotic approximation. In this case, the error is larger and, up to small oscillations, $\frac{t}{\ln t} e_r(t) := \frac{t}{\ln t} \sup_{x \in [10,70]}|U(x,t) - u_{\lead,r}(x,t)|$ is a growing function of $t$. 

This shows that the numerical scheme is able to capture the left-moving wave, and illustrates the conclusion drawn in \cite{CLmain} that the ``bad'' Boussinesq equation has a region where the nonlinear interaction between right- and left-moving waves makes a quantitative difference even asymptotically.

\begin{figure}
\bigskip\begin{center}
\hspace{-.4cm}
\begin{overpic}[width=.46\textwidth]{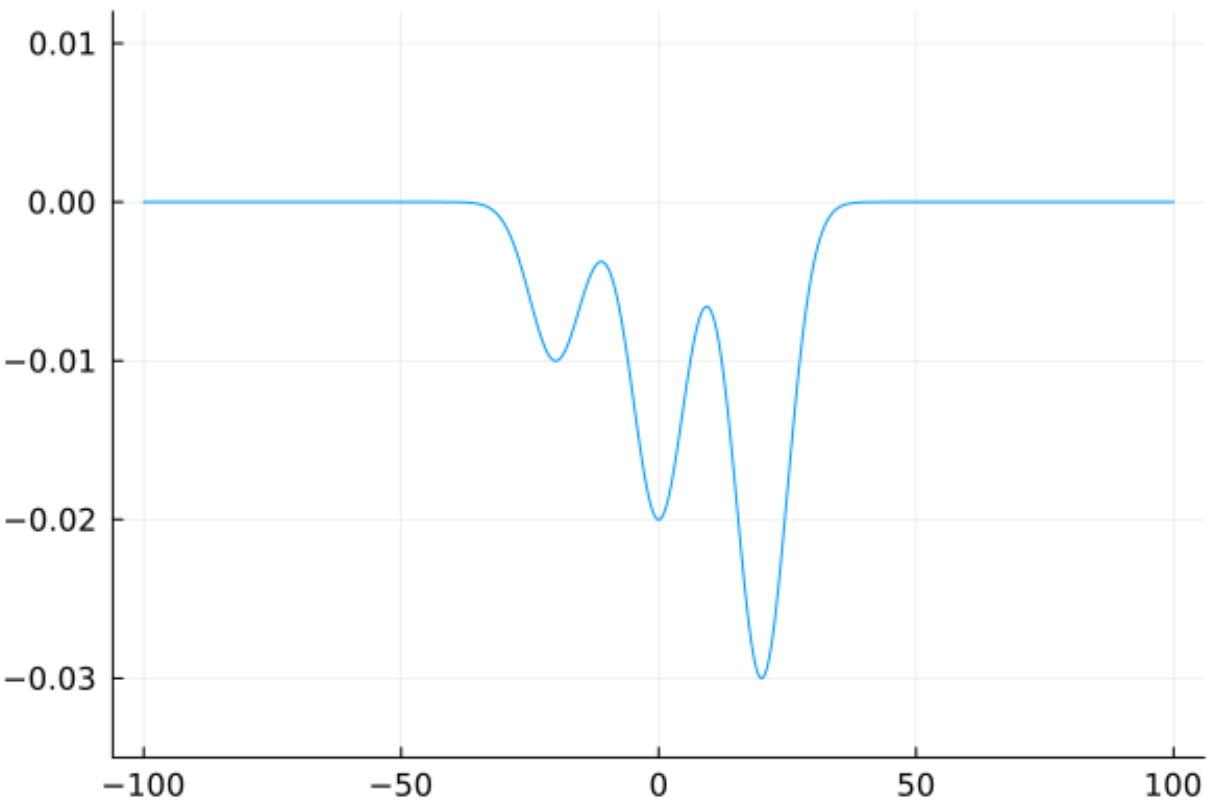}
      \put(8,68){\footnotesize $U$}
      \put(101.6,2.6){\footnotesize $x$}
       \put(49,54){\footnotesize $t = 0$}
    \end{overpic}
    \hspace{0.35cm}
\begin{overpic}[width=.46\textwidth]{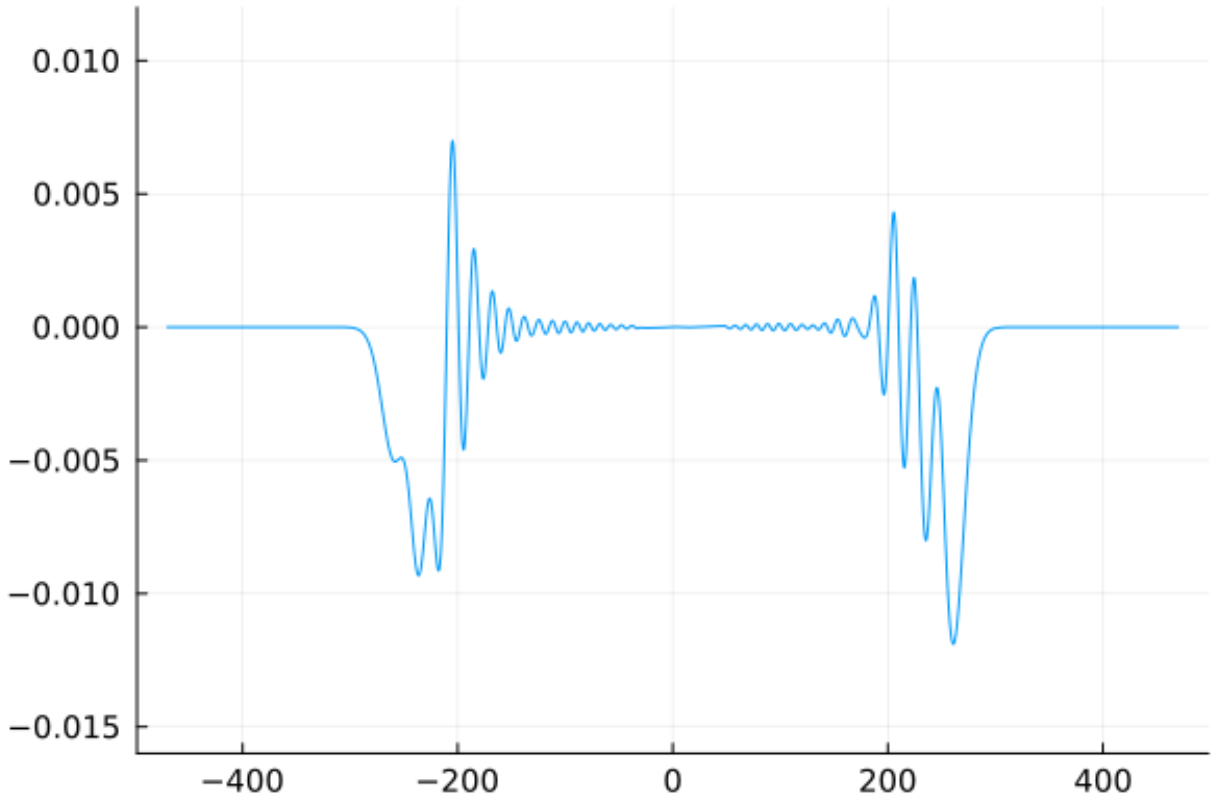}
      \put(8,68){\footnotesize $U$}
      \put(101.6,2.8){\footnotesize $x$}
       \put(47,48){\footnotesize $t = 250$}
    \end{overpic}
    \\
    \vspace{0.6cm}
    \hspace{-0.5cm}
\begin{overpic}[width=.46\textwidth]{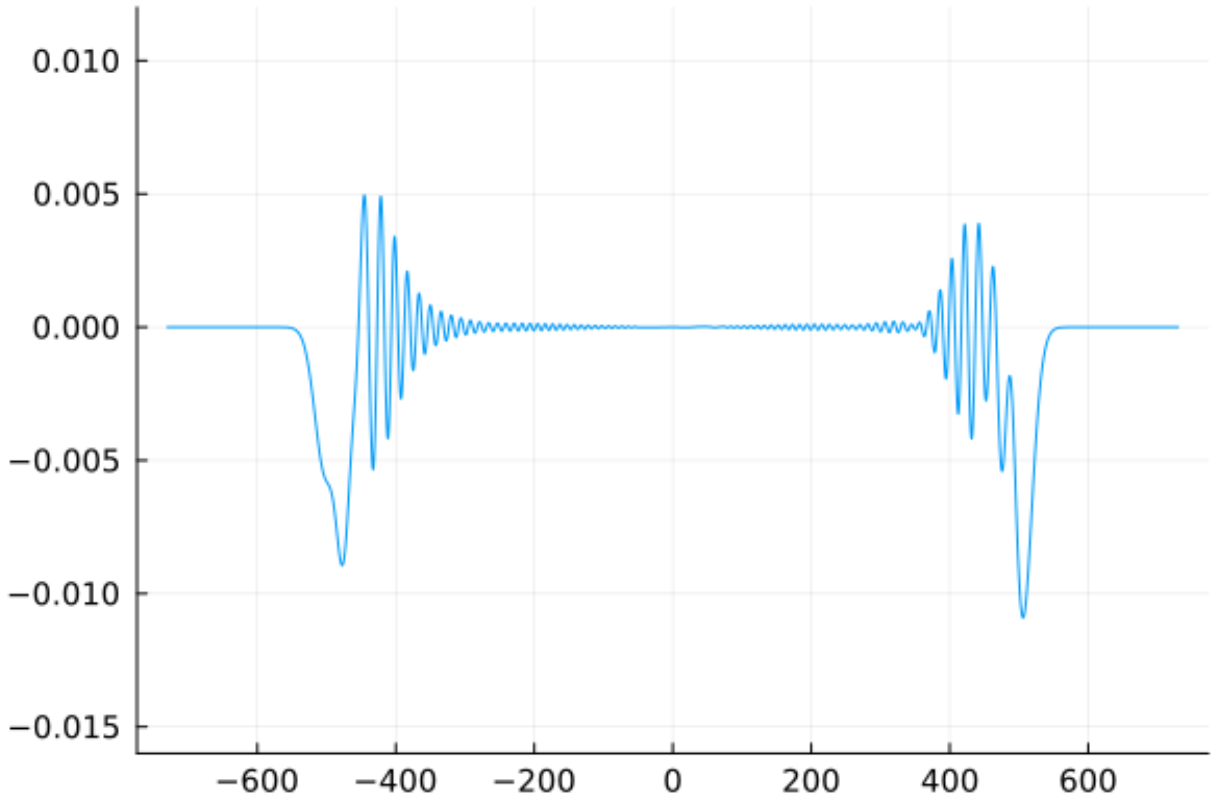}
      \put(10,68){\footnotesize $U$}
      \put(101.6,2.8){\footnotesize $x$}
       \put(47,48){\footnotesize $t = 500$}
    \end{overpic}
        \hspace{0.54cm}
\begin{overpic}[width=.46\textwidth]{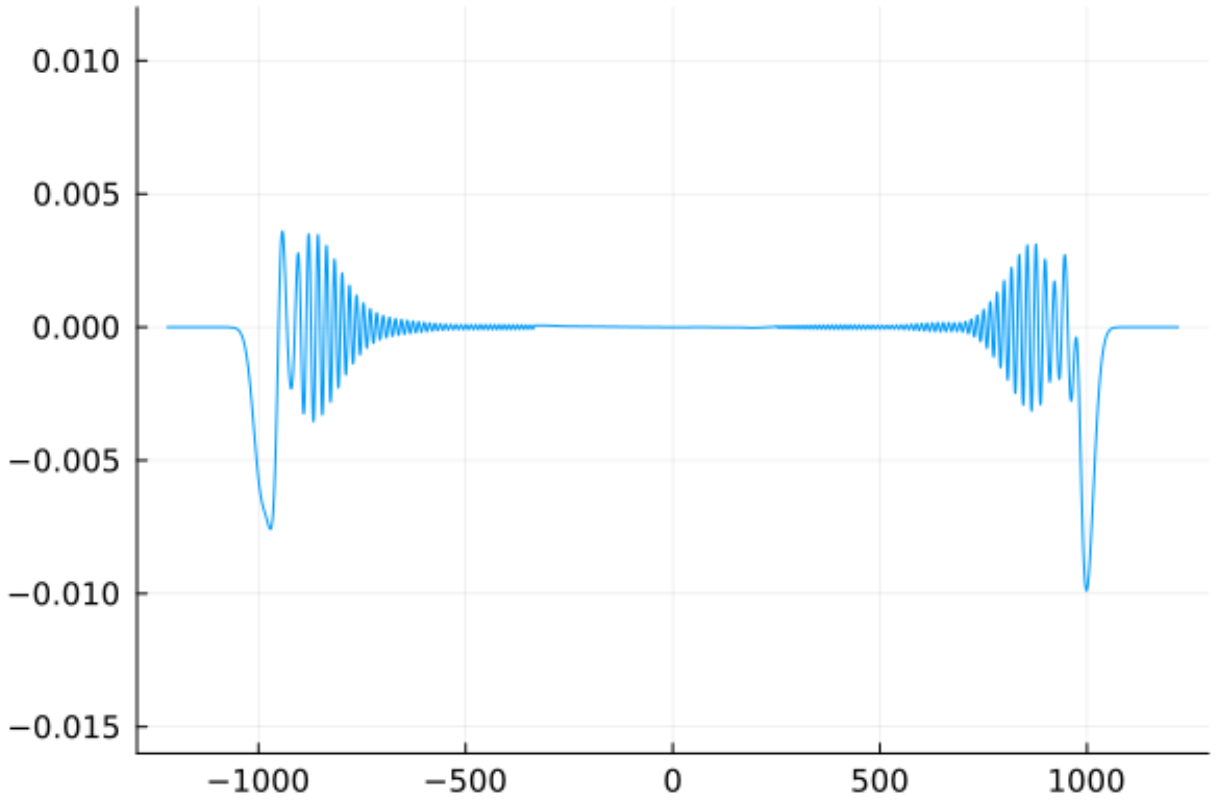}
      \put(10,68){\footnotesize $U$}
      \put(101.6,2.8){\footnotesize $x$}
       \put(47,48){\footnotesize $t = 1000$}
    \end{overpic}
     \begin{figuretext}\label{threehumpsfig1}
       The approximate solution $U(x,t)$ of (\ref{badboussinesq}) obtained by the scheme of Section \ref{schemesec} for the initial data (\ref{threehumps}) at times $t = 0$, $t=250$, $t=500$, and $t=1000$. 
\end{figuretext}
     \end{center}
\end{figure}

\begin{figure}
\bigskip\begin{center}
\hspace{-.50cm}
\begin{overpic}[width=.46\textwidth]{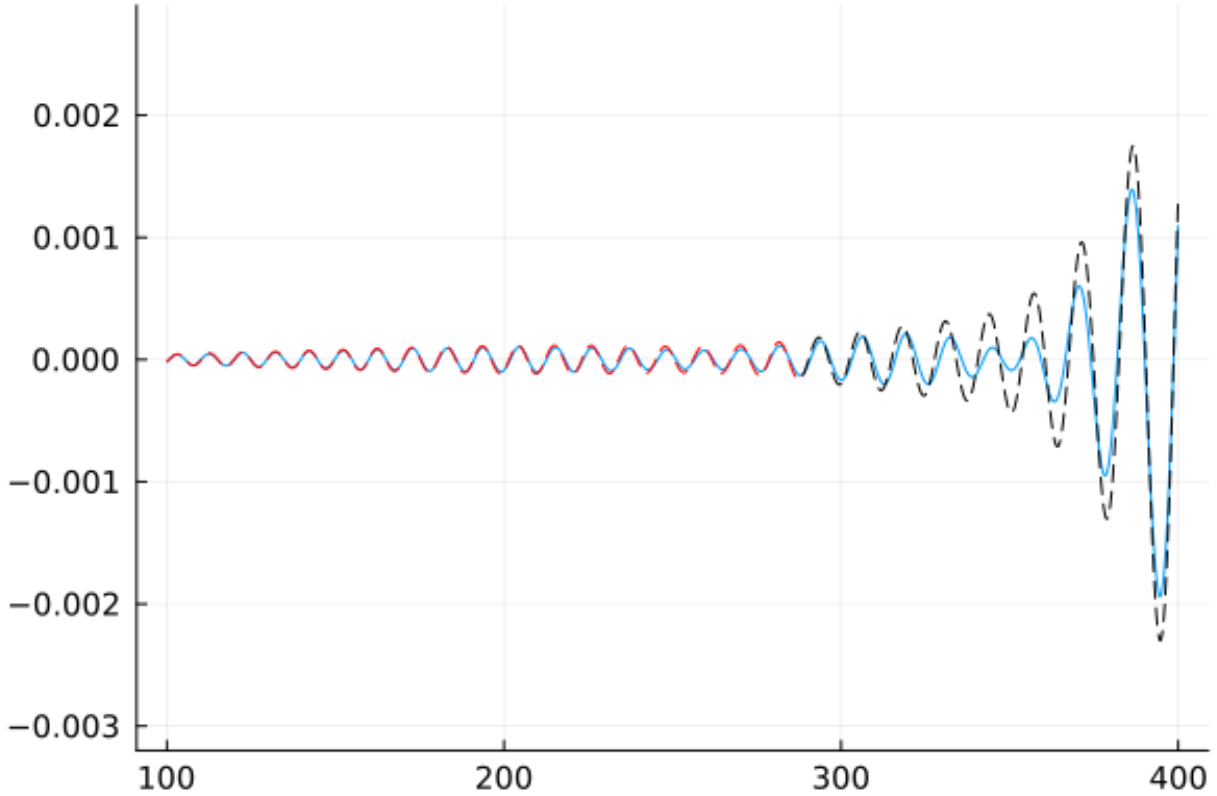}
      \put(3,68){\footnotesize $u_{\lead}, U$}
      \put(101.6,2.2){\footnotesize $x$}
       \put(47,48){\footnotesize $t = 500$}
    \end{overpic}
    \hspace{0.4cm}
\begin{overpic}[width=.46\textwidth]{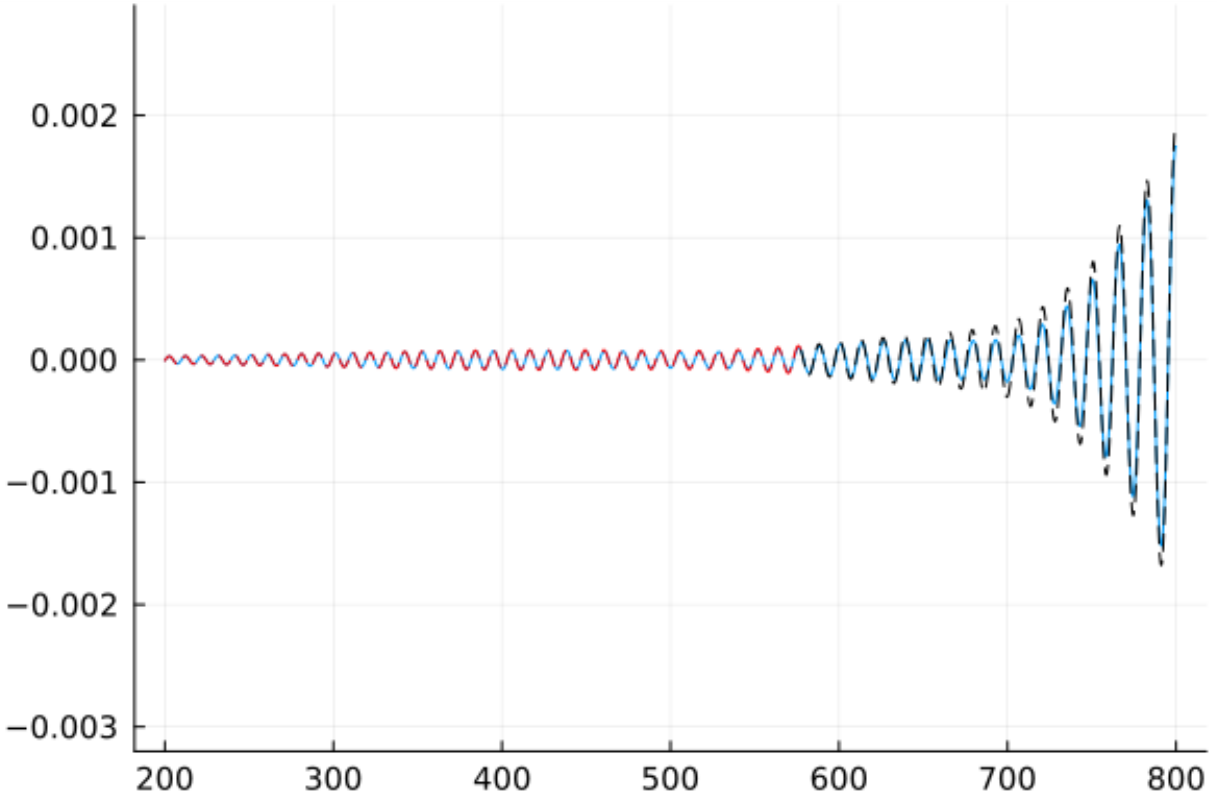}
      \put(3,68){\footnotesize $u_{\lead}, U$}
      \put(101.6,2.2){\footnotesize $t$}
       \put(47,48){\footnotesize $t = 1000$}
    \end{overpic}
	\\        \vspace{0.7cm}
    \hspace{-0.20cm}
\begin{overpic}[width=.46\textwidth]{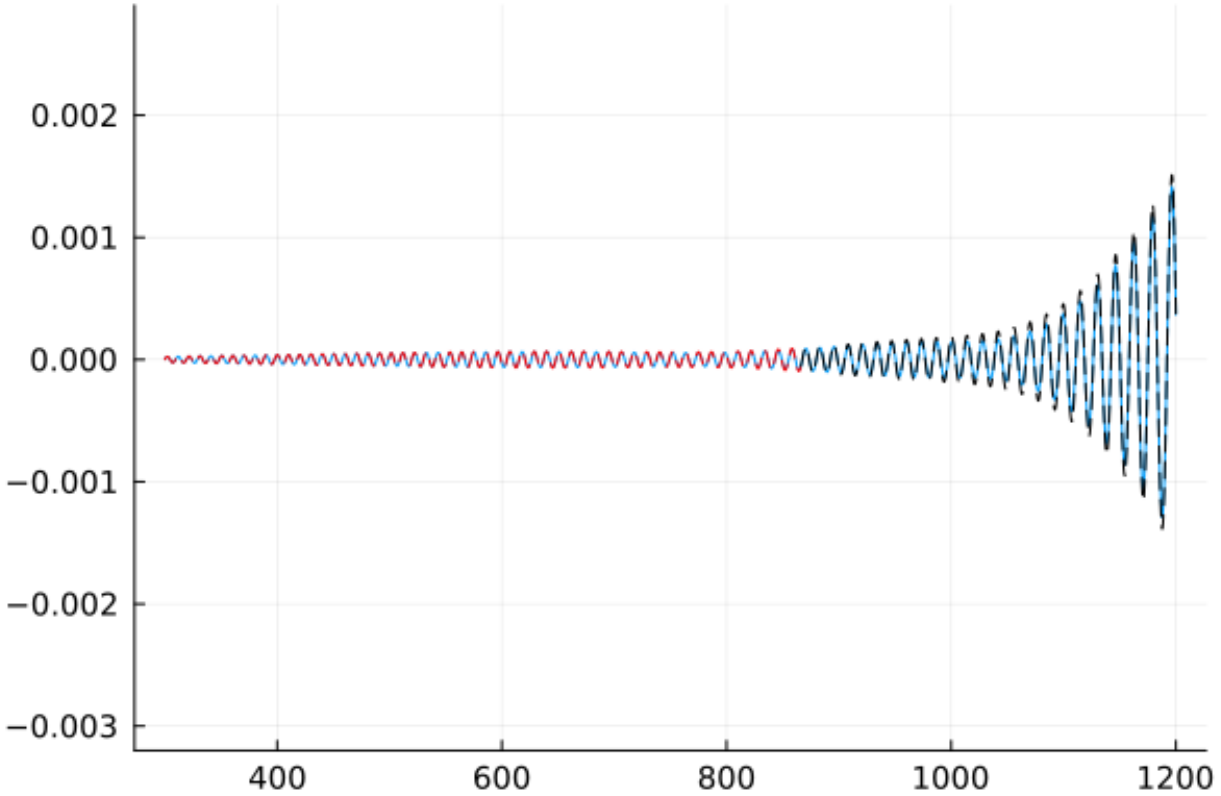}
      \put(3,68){\footnotesize $u_{\lead}, U$}
      \put(101.6,2.2){\footnotesize $x$}
       \put(47,48){\footnotesize $t = 1500$}
    \end{overpic}
    \hspace{0.6cm}
\begin{overpic}[width=.46\textwidth]{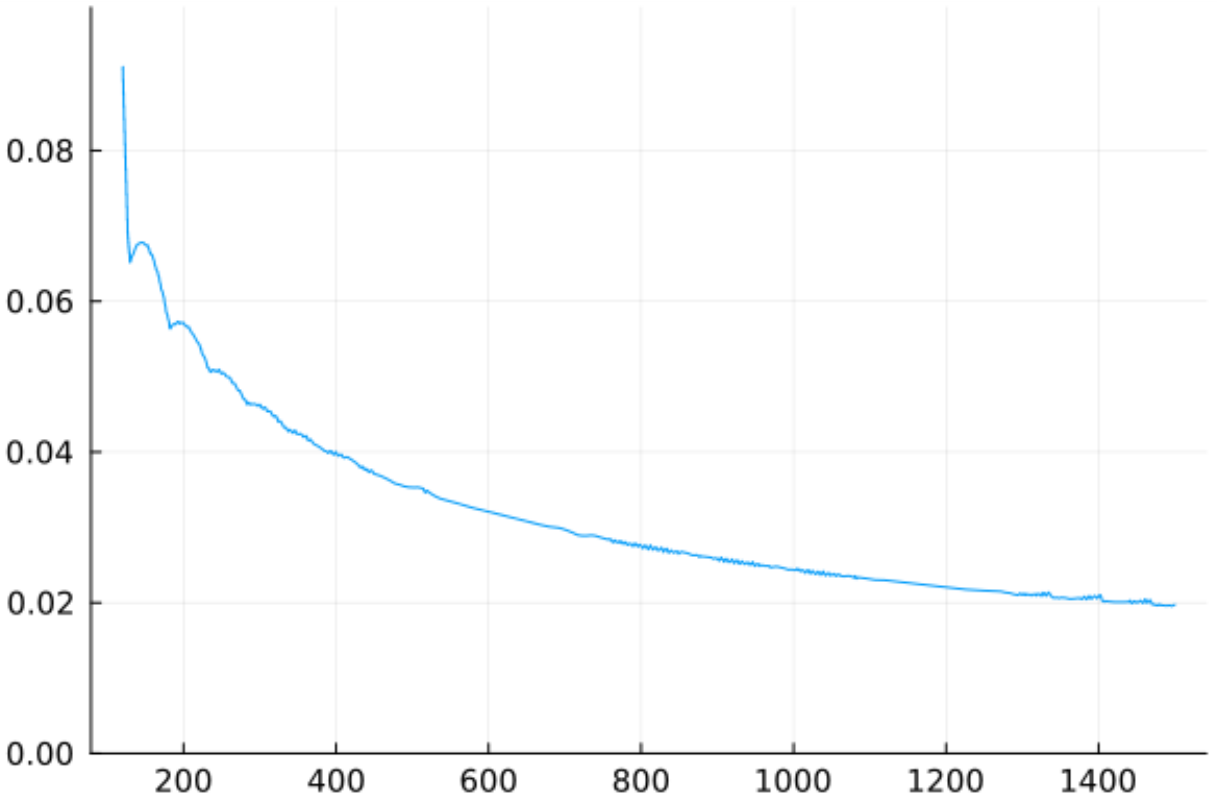}
      \put(0,69){\footnotesize $\frac{t}{\ln t}e(t)$}
      \put(101.6,2.2){\footnotesize $t$}
    \end{overpic}
     \begin{figuretext}\label{threehumpsfig2}
       Topleft, topright, and bottomleft: The leading asymptotic term $u_{\lead}(x,t)$ (dashed black in Sector IV and dashed red in Sector V) together with the numerical approximation $U(x,t)$ (solid blue) for the initial data (\ref{threehumps}) obtained with the scheme of Section \ref{schemesec} at $t = 500$, $t = 1000$, and $t = 1500$, respectively. 
       
       Bottomright: The graph of $t \mapsto \frac{t}{\ln t} e(t)$ where $e(t) = \sup_{x \in [0.2t, 0.8t]} |U(x,t) - u_{\lead}(x,t)|$.
\end{figuretext}
     \end{center}
\end{figure}

\subsection{A sum of three Gaussians}\label{threehumpssubsec}
We next consider the initial data 
\begin{align}\label{threehumps}
u(x,0) = -3a e^{-c(x-b)^2}-2ae^{-c x^2} - a e^{-c(x+b)^2}, \qquad u_t(x,0) = 0,
\end{align}
with $a = 0.01$, $b = 20$, and $c=0.02$.
Applying the scheme of Section \ref{schemesec} to this initial data, we obtain the approximate solution $U(x,t)$ shown in Figure \ref{threehumpsfig1}. 
Figure \ref{threehumpsfig2} shows $u_{\lead}(x,t)$ and the numerical approximation $U(x,t)$ for the initial data (\ref{threehumps}) at three different times, as well as the graph of $t \mapsto \frac{t}{\ln t} e(t)$ where $e(t) = \sup_{x \in [0.2t, 0.8t]} |U(x,t) - u_{\lead}(x,t)|$. The function $\frac{t}{\ln t} e(t)$ is roughly twice as large for (\ref{threehumps}) as it was in the case of the single Gaussian (\ref{Gaussian}) (cf. Figure \ref{onehumpfig2}), but it still exhibits no tendency to increase with time, in consistency with (\ref{uasymptotics}).

\subsection{A solution supporting a soliton}
\begin{figure}
\bigskip\begin{center}
\hspace{-.4cm}
\begin{overpic}[width=.46\textwidth]{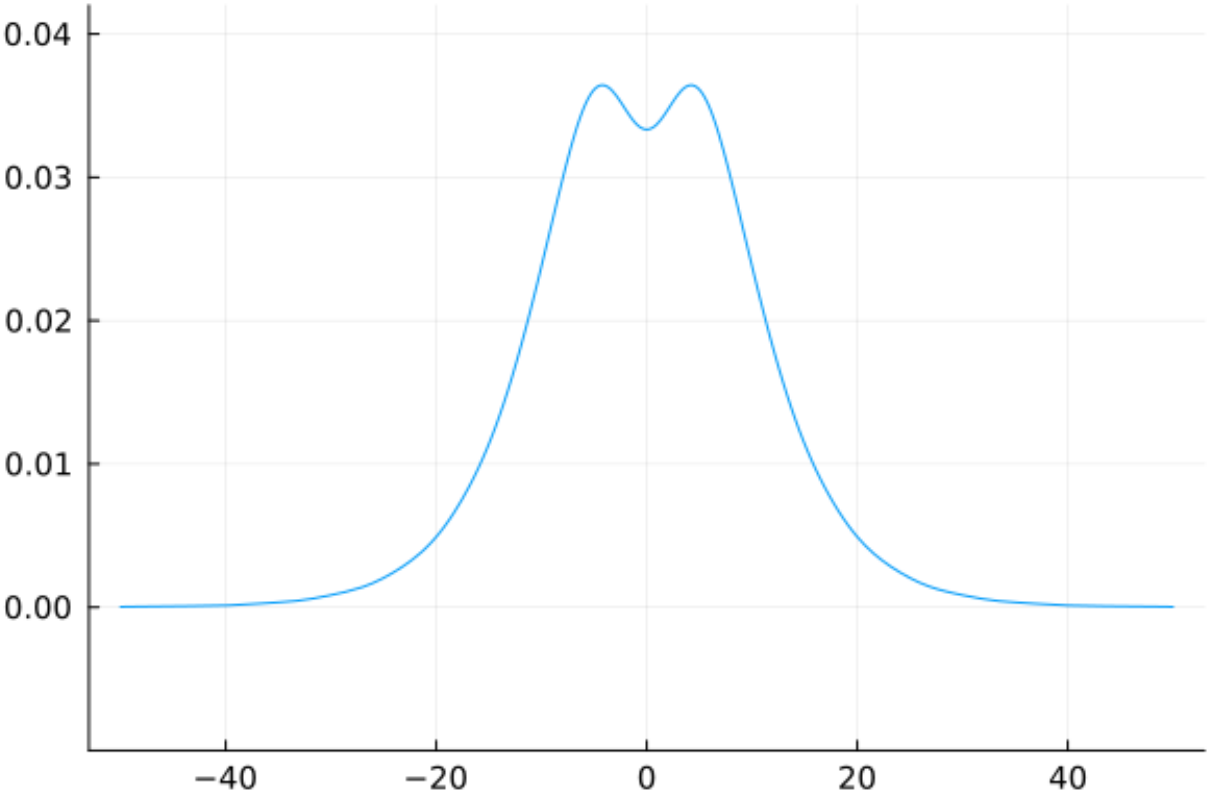}
      \put(6,68){\footnotesize $U$}
      \put(101.6,2.6){\footnotesize $x$}
       \put(48,43){\footnotesize $t = 0$}
    \end{overpic}
    \hspace{0.35cm}
\begin{overpic}[width=.46\textwidth]{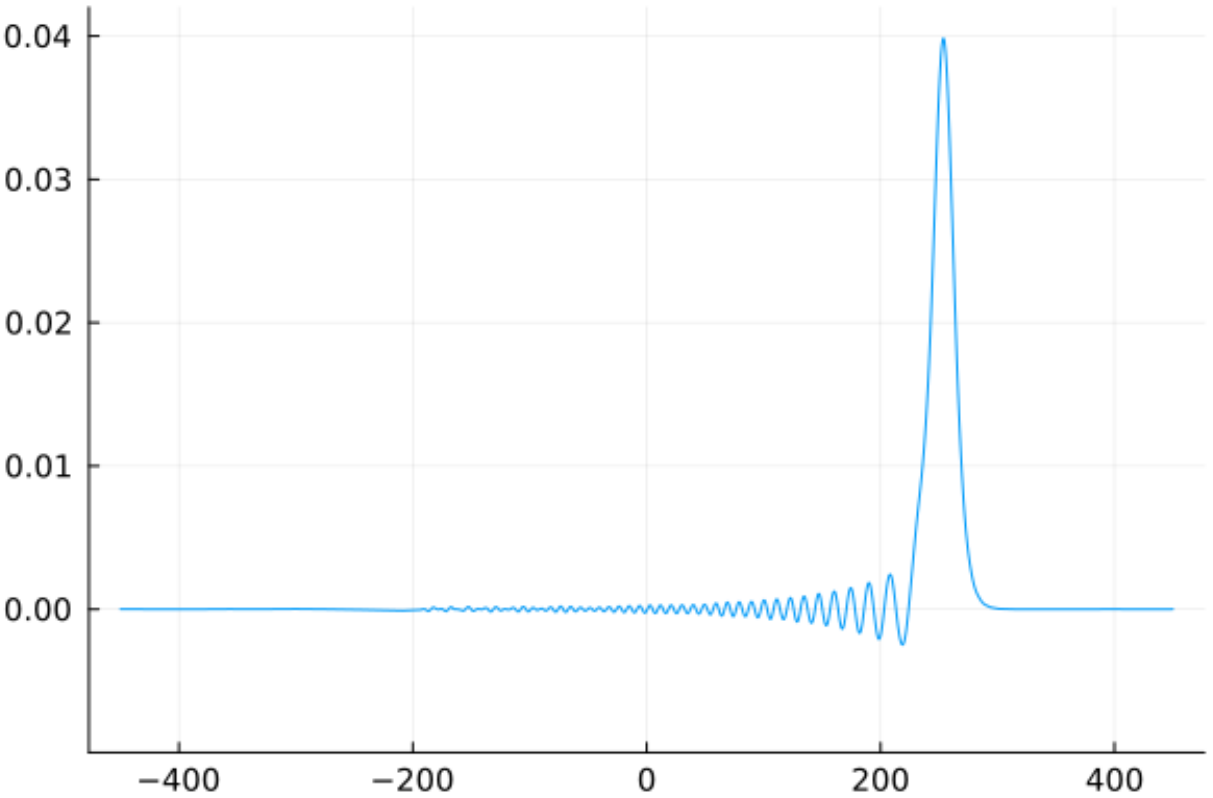}
      \put(6,68){\footnotesize $U$}
      \put(101.6,2.8){\footnotesize $x$}
       \put(47,43){\footnotesize $t = 250$}
    \end{overpic}
    \\
    \vspace{0.6cm}
    \hspace{-0.5cm}
\begin{overpic}[width=.46\textwidth]{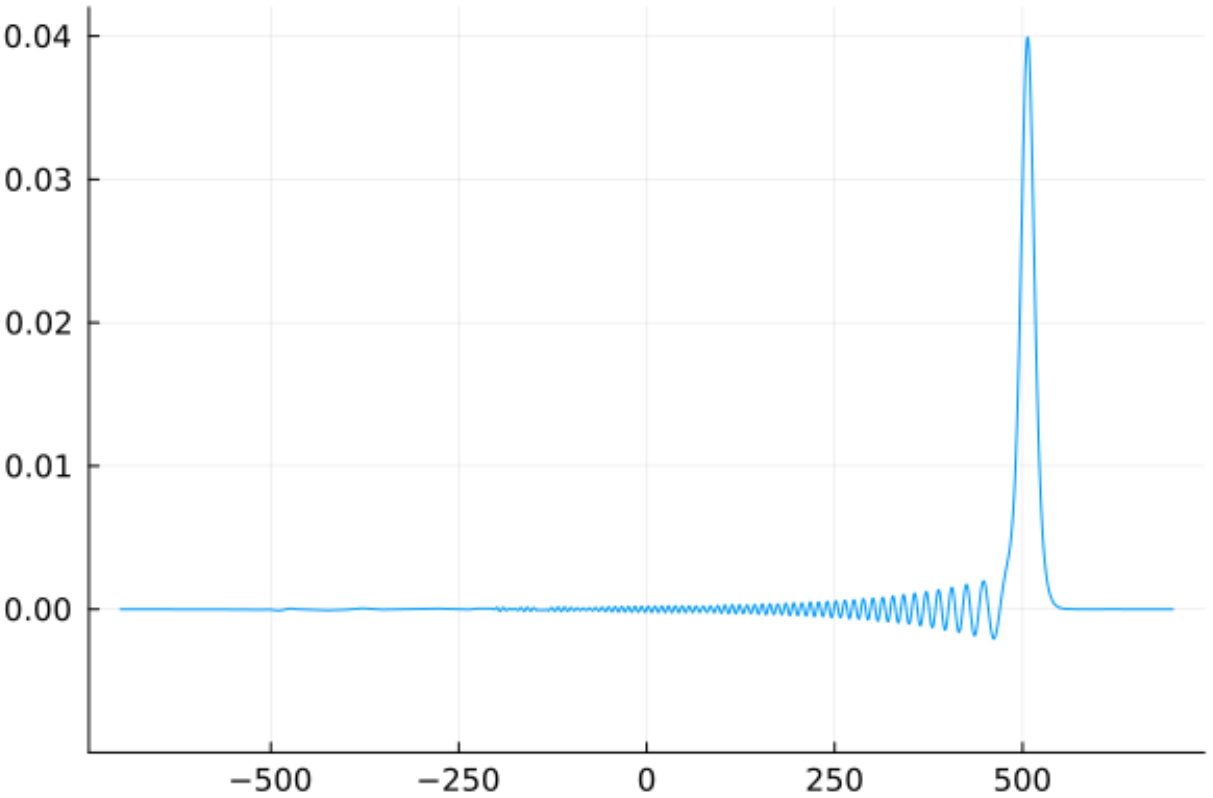}
      \put(6,68){\footnotesize $U$}
      \put(101.6,2.8){\footnotesize $x$}
       \put(47,43){\footnotesize $t = 500$}
    \end{overpic}
        \hspace{0.54cm}
\begin{overpic}[width=.46\textwidth]{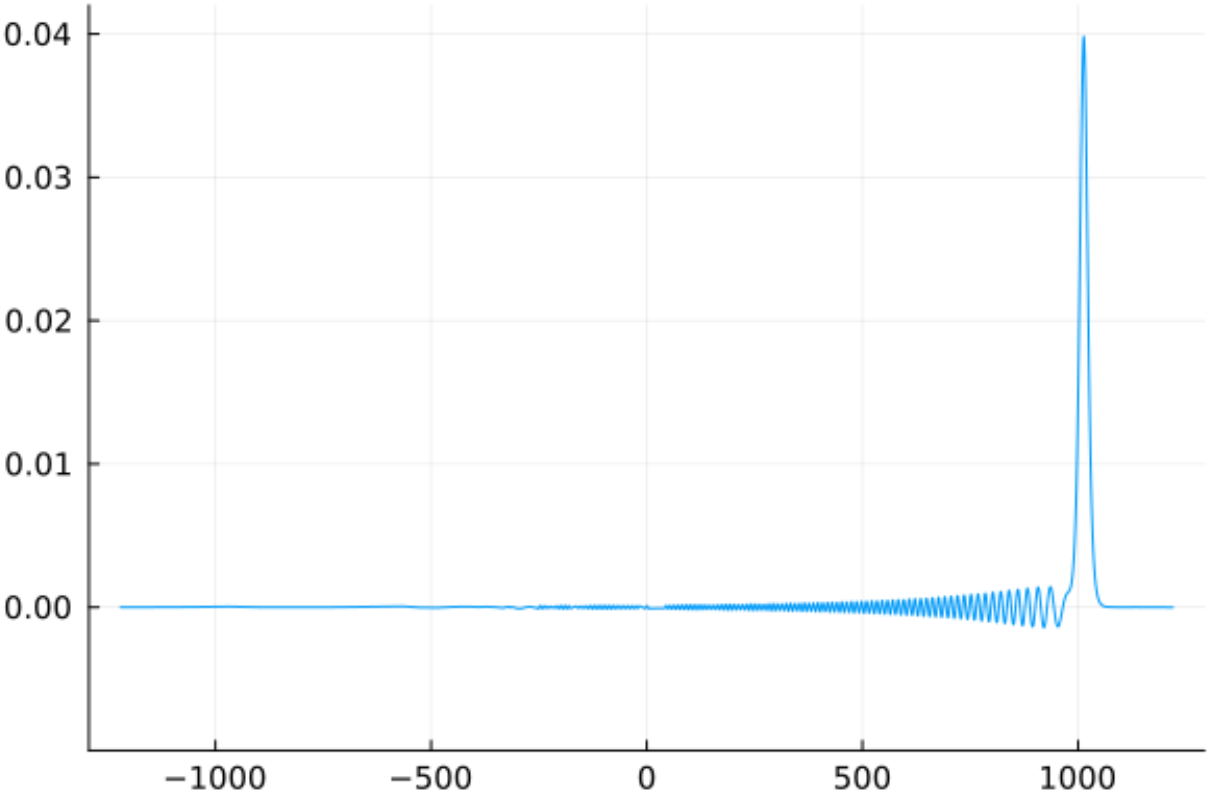}
      \put(6,68){\footnotesize $U$}
      \put(101.6,2.8){\footnotesize $x$}
       \put(47,43){\footnotesize $t = 1000$}
    \end{overpic}
     \begin{figuretext}\label{solutionwithsolitonfig1}
       The approximate solution $U(x,t)$ of (\ref{badboussinesq}) obtained by the scheme of Section \ref{schemesec} for the initial data (\ref{solutionwithsolitoninitialdata}) at times $t = 0$, $t=250$, $t=500$, and $t=1000$. The right-moving hump with amplitude approximately $0.04$ is the asymptotic soliton generated by the zero $k_0 \approx 1.1755$ of $s_{11}(k)$. The asymptotic soliton is followed by a decaying dispersive tail. 
\end{figuretext}
     \end{center}
\end{figure}

The asymptotic formula (\ref{uasymptotics}) is valid for solitonless solutions, i.e., for solutions such that the spectral function $s_{11}(k)$ defined in (\ref{sdef}) has no zeros. 
%for $k \in (\bar{D}_2 \cup \partial \D) \setminus (\{e^{\frac{\pi i(j-1)}{3}}\}_{j=1}^{6} \cup \{0\})$, where $D_2 = \{k \in \C\,|\, \re l_1 < \re l_3 < \re l_2\}$.
If $s_{11}(k)$ has a zero at a point $k_0 \in (-1,0) \cup (1,\infty)$, then the solution $u(x,t)$ develops a soliton for large times; zeros of $s_{11}(k)$ in $(-1,0)$ generate left-moving solitons, whereas zeros in $(1,+\infty)$ generate right-moving solitons, see \cite{CLscatteringsolitons}. The solutions considered in Sections \ref{onehumpsubsec} and \ref{threehumpssubsec} are all solitonless. In this section, we consider a solution that supports a soliton. More precisely, we let the initial data $u(x,0)$ be the sum of the initial data of the exact one-soliton considered in Section \ref{solitonsec}, $A \sech^2(\sqrt{A/6}x)$, and a Gaussian perturbation, $-A e^{-A x^2}/3$:
\begin{align}\label{solutionwithsolitoninitialdata}
u(x,0) = A \sech^2\bigg(\sqrt{\frac{A}{6}}x\bigg) - \frac{A}{3} e^{-A x^2}, \qquad u_t(x,0) = -c u_x(x,0),
\end{align}
where $A = 0.05$ and $c = \sqrt{1 + \frac{2A}{3}}$. %If the Gaussian perturbation were absent, then the corresponding solution $u(x,t)$ of the initial-value problem for (\ref{badboussinesq}) would reduce to the exact one-soliton $A \sech^2(\sqrt{A/6}(x - ct))$ studied in Section \ref{solitonsec}. 
As discussed below, the corresponding solution $u(x,t)$ of the initial-value problem for (\ref{badboussinesq}) asymptotes for large times to a right-moving one-soliton generated by the zero $k_0 \approx 1.1755 >1$ of $s_{11}(k)$.
The scheme of Section \ref{schemesec} applied to the initial data (\ref{solutionwithsolitoninitialdata}) yields the approximate solution $U(x,t)$ shown in Figure \ref{solutionwithsolitonfig1}. 

\subsubsection{Comparison with the asymptotic approximation in Sectors IV and V}
Let us first compare the approximate solution $U(x,t)$ with the leading asymptotic approximation in Sectors IV and V; as in Section \ref{onehumpsubsec}, these sectors are given by $\frac{1}{\sqrt{3}} < \zeta < 1$ and $0 < \zeta < \frac{1}{\sqrt{3}}$, respectively, where $\zeta = x/t$. In the presence of a soliton, the asymptotic formula (\ref{uasymptotics}) for $u(x,t)$ is modified as follows. Let $k_2 = k_2(\zeta)$ be as in (\ref{k2def}), let
$$k_{4} = k_{4}(\zeta) = \frac{1}{4}\bigg( \zeta + \sqrt{8+\zeta^{2}} - \sqrt{2}\sqrt{-4+\zeta^{2}+\zeta\sqrt{8+\zeta^{2}}} \bigg),$$
and define the rational function $\mathcal{P}(k)$ by
$$\mathcal{P}(k)=\frac{(k-\omega^2 k_0)(k-\omega k_0^{-1})}{(k- \omega k_0)(k-\omega^2k_0^{-1})},$$ 
where $\omega = e^{\frac{2\pi i}{3}}$. The solution $u(x,t)$ with initial data (\ref{solutionwithsolitoninitialdata}) obeys (\ref{uasymptotics}) with the asymptotic term $u_{\lead}(x,t)$ modified by introducing phase shifts $ \arg \tfrac{\mathcal{P}(\omega k_{4})}{\mathcal{P}(\omega^{2} k_{4})} $ and $\arg \tfrac{\mathcal{P}(\omega^{2} k_{2})}{\mathcal{P}(\omega k_{2})}$ in the cosine factors as follows \cite{CLsectorIV}: 
\begin{align} \nonumber
u_{\lead}(x,t) = 
&\; \frac{A_{1}(\zeta)}{\sqrt{t}}   \cos\Big(\alpha_{1}(\zeta,t) + \arg \tfrac{\mathcal{P}(\omega k_{4})}{\mathcal{P}(\omega^{2} k_{4})} \Big) 
	\\ \label{uleaddefsol1}
& + \frac{A_{2}(\zeta)}{\sqrt{t}} \cos\Big( \alpha_{2}(\zeta,t) + \arg \tfrac{\mathcal{P}(\omega^{2} k_{2})}{\mathcal{P}(\omega k_{2})}\Big) \qquad  \text{in Sector IV},
\end{align}
and
\begin{align} \nonumber
u_{\lead}(x,t) =
&\; \frac{\tilde{A}_{1}(\zeta)}{\sqrt{t}}     \cos\Big(\tilde{\alpha}_{1}(\zeta,t) + \arg \tfrac{\mathcal{P}(\omega k_{4})}{\mathcal{P}(\omega^{2} k_{4})} \Big)
	\\ \label{uleaddefsol2}
&  +\frac{A_{2}(\zeta)}{\sqrt{t}}  \cos\Big( \tilde{\alpha}_{2}(\zeta,t) + \arg \tfrac{\mathcal{P}(\omega^{2} k_{2})}{\mathcal{P}(\omega k_{2})}\Big) \qquad  \text{in Sector V},
 \end{align}
where $A_1, A_2, \alpha_1, \alpha_2, \tilde{A}_1, \tilde{\alpha}_1, \tilde{\alpha}_2$ are as in (\ref{uleaddef}).
%The phase shifts describe leading order soliton-radiation interaction in sector IV and V. 
When the solution supports no soliton, $\mathcal{P}(k)$ is identically 1. In that case, the phase shifts are absent and (\ref{uleaddefsol1})--(\ref{uleaddefsol2}) reduce to (\ref{uleaddef}). 

\begin{figure}
\bigskip\begin{center}
\hspace{-.50cm}
\begin{overpic}[width=.46\textwidth]{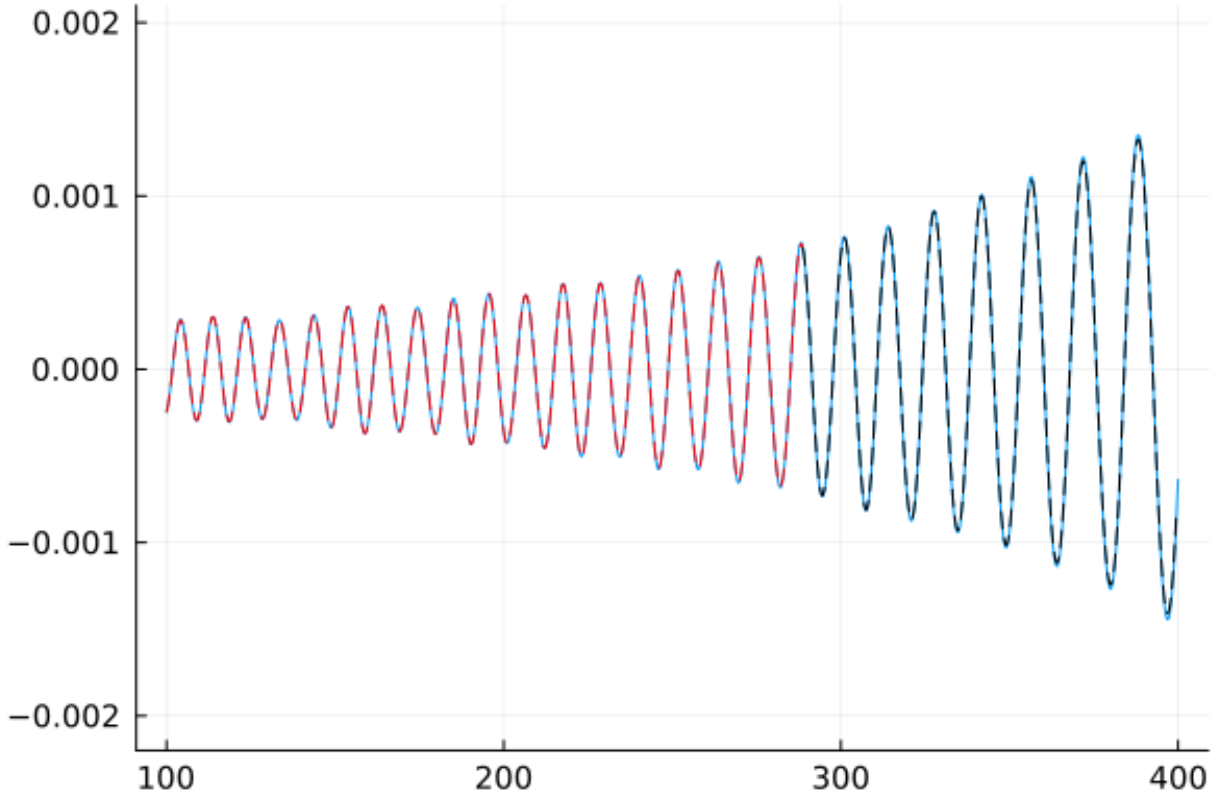}
      \put(3,68){\footnotesize $u_{\lead}, U$}
      \put(101.6,2.2){\footnotesize $x$}
       \put(47,55){\footnotesize $t = 500$}
    \end{overpic}
        \hspace{0.4cm}
\begin{overpic}[width=.46\textwidth]{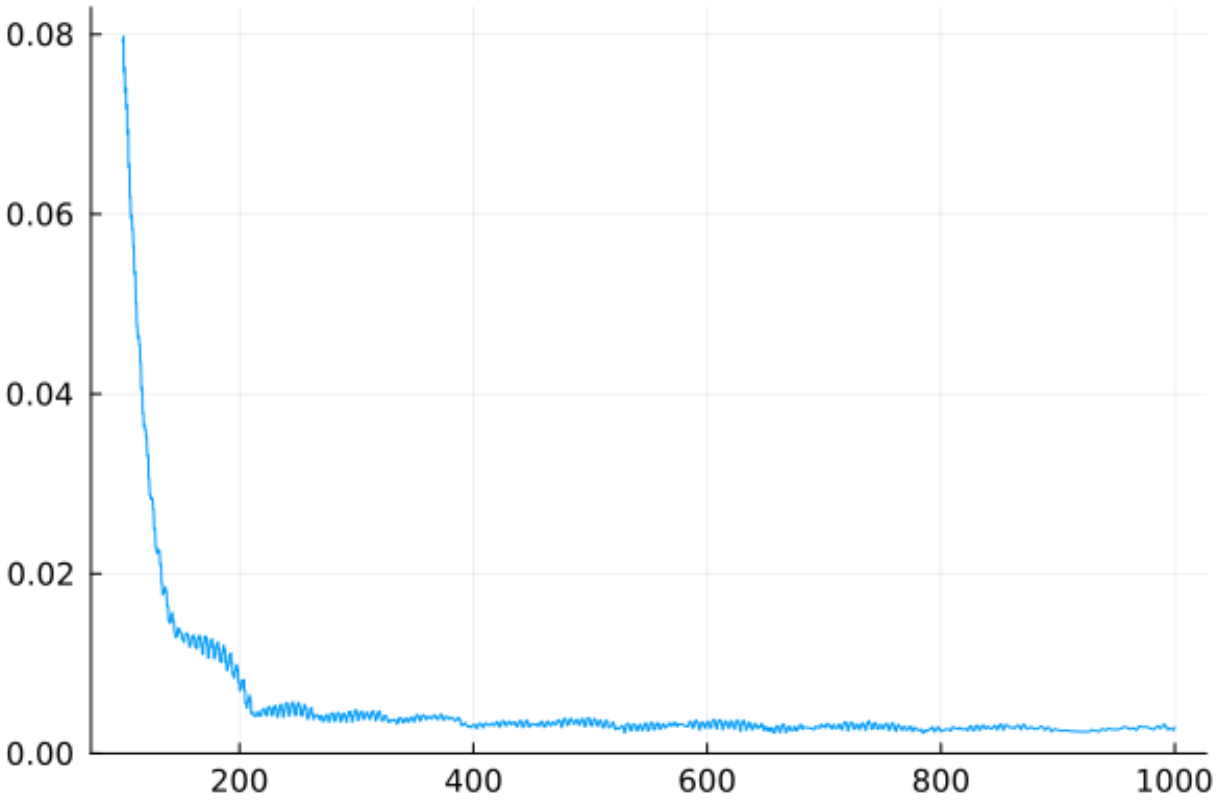}
      \put(0,69){\footnotesize $\frac{t}{\ln t}e(t)$}
      \put(101.6,2.2){\footnotesize $t$}
    \end{overpic}
     \begin{figuretext}\label{solutionwithsolitonfig2}
       Left: The asymptotic approximation $u_{\lead}(x,t)$ for $x \in [0.2t, 0.8t]$ and $t = 500$ (dashed black in Sector IV and dashed red in Sector V) together with the numerical approximation $U(x,t)$ (solid blue) for the initial data (\ref{solutionwithsolitoninitialdata}) obtained with the scheme of Section \ref{schemesec}. 
       
       Right: The graph of $t \mapsto \frac{t}{\ln t} e(t)$ where $e(t) = \sup_{x \in [0.2t, 0.8t]} |U(x,t) - u_{\lead}(x,t)|$.
\end{figuretext}
     \end{center}
\end{figure}

Figure \ref{solutionwithsolitonfig2} (left) shows $u_{\lead}(x,t)$ and the numerical approximation $U(x,t)$ of the solution $u(x,t)$ for $x \in [0.2t, 0.8t]$ at $t = 500$.
The theoretical results of \cite{CLsectorV, CLsectorIV} show that $u_{\lead}(x,t)$ approximates $u(x,t)$ in Sectors IV  and  V with an error of order $O(\frac{\ln t}{t})$ for large $t$. In Figure \ref{solutionwithsolitonfig2} (right) we plot the graph of the $L^\infty$-error $e(t) = \sup_{x \in [0.2t, 0.8t]} |U(x,t) - u_{\lead}(x,t)|$ multiplied by $\frac{t}{\ln t}$. If $U(x,t)$ is a good approximation of $u(x,t)$, then $\frac{t}{\ln t} e(t)$ should be of $O(1)$ for large $t$; Figure \ref{solutionwithsolitonfig2} confirms this expectation.

\begin{figure}
\bigskip
\begin{center}
\hspace{-.4cm}
\begin{overpic}[width=.46\textwidth]{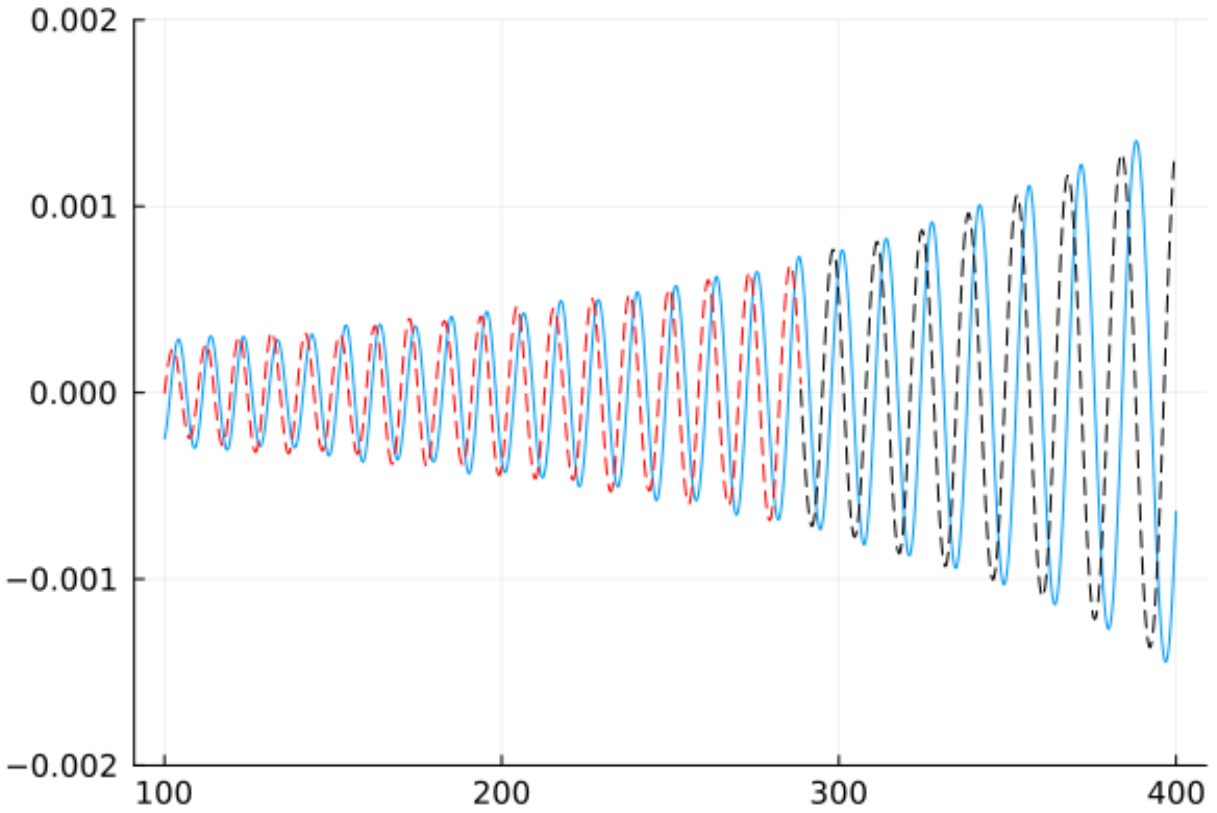}
      \put(4,69){\footnotesize $u_{\lead, \noshift}, U$}
      \put(101.6,3){\footnotesize $x$}
       \put(47,55){\footnotesize $t = 500$}
    \end{overpic}
    \hspace{0.5cm}
\begin{overpic}[width=.46\textwidth]{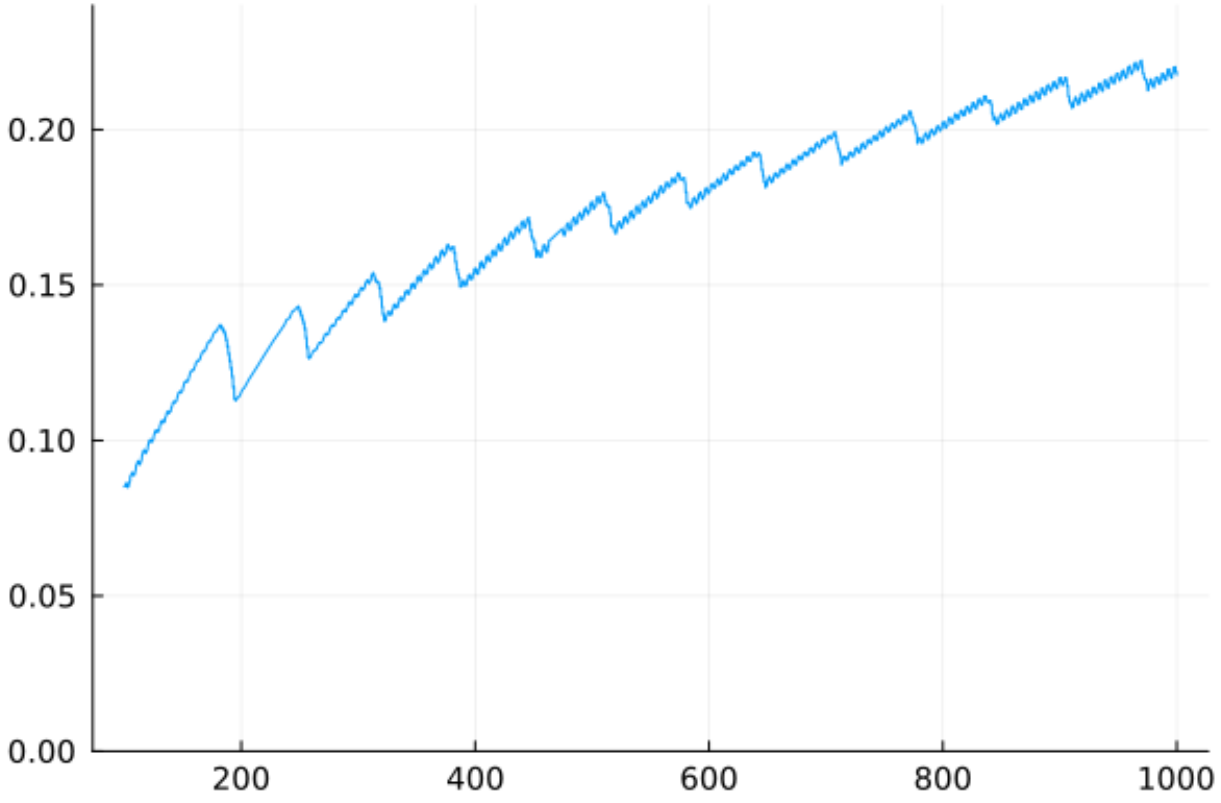}
      \put(1,69){\footnotesize $\frac{t}{\ln t}e(t)$}
      \put(101.6,3){\footnotesize $t$}
    \end{overpic}
     \begin{figuretext}[Large error without phase shifts]\label{solitonnophase}
       Left: The asymptotic approximation $u_{\lead, \noshift}(x,t)$ without phase shifts (dashed black in Sector V and red in sector IV) together with the numerical approximation $U(x,t)$ (solid blue) obtained with the scheme of Section \ref{schemesec} at time $t = 500$.
       
       Right: The graph of $t \mapsto \frac{t}{\ln t} e(t) = \frac{t}{\ln t} \sup_{x \in [0.2t,0.8t]}|U(x,t) - u_{\lead,\noshift}(x,t)|$.
\end{figuretext}
     \end{center}
\end{figure}

To illustrate the importance of the phase shifts, we let $u_{\lead, \noshift}(x,t)$ denote $u_{\lead}(x,t)$ without the phase shifts. Figure \ref{solitonnophase} (left) shows $u_{\lead, \noshift}(x,t)$ for $t=500$ and Figure \ref{solitonnophase} (right) displays the corresponding $L^\infty$-error. The phase shifts are clearly of crucial importance.

\subsubsection{Comparison with the asymptotic approximation in Sector II}
We next compare the approximate solution $U(x,t)$ with the leading asymptotic approximation in the asymptotic sector $1 + \epsilon \leq \zeta = x/t \leq M$, where $\epsilon > 0$ and $M > 0$ are arbitrary constants. This sector was referred to as Sector II in \cite{CLmain}.
The solution $u(x,t)$ with initial data (\ref{solutionwithsolitoninitialdata}) obeys the following asymptotic formula in Sector II \cite{CLsectorII}:
\begin{align}\label{uleadingsectorII}
u(x,t) = u_{\sol}(x,t) + O(t^{-1/2}) \qquad \text{as $t \to +\infty$}
\end{align}
uniformly for $x/t \in [1+\epsilon, M]$, where 
\begin{equation} \label{usoldef} 
  u_\sol(x,t)= A_0 \sech^2\bigg(\sqrt{\frac{A_0}{6}}(x-c_0 t)-\ln f_{k_0}(\zeta)\bigg).
\end{equation} 
The real constants $A_0$ and $c_0$ in (\ref{usoldef}) are defined by
\begin{align*}
& A_0 :=\frac{3}{8}(k_0-k_0^{-1})^2 \approx 0.03955, \qquad c_0 := \sqrt{1 + \frac{2A_0}{3}} =\frac{k_0+k_0^{-1}}{2} \approx 1.0131, 
\end{align*}
where $k_0 \approx 1.1755 >1$ is the zero of $s_{11}(k)$.
The definition of the real-valued function $\ln f_{k_0}(\zeta)$ in (\ref{usoldef}) is given in the appendix where we also describe how equation (\ref{uleadingsectorII}) follows from the results of \cite{CLscatteringsolitons, CLsectorII}.

As $t \to \infty$, the leading term $u_\sol(x,t)$ in (\ref{uleadingsectorII}) asymptotes to the exact one-soliton of amplitude $A_0$ and speed $c_0$ given by
$$A_0 \sech^2\bigg(\sqrt{\frac{A_0}{6}}(x-c_0 t - x_0)\bigg) \qquad \text{with $x_0 := \sqrt{\frac{6}{A_0}} \ln f_{k_0}(c_0)$}.$$ %x_0 \approx -0.11323
Indeed, for large times the $\sech^2$-function in (\ref{usoldef}) is exponentially suppressed except in a narrow region around the line  $x/t = c_0$ whose width tends to zero as $t \to +\infty$. 
%In our example, $\Delta(x/t,k_0)$ is close to 1, so $x_0\approx 0$.

\begin{figure}[h!]
\bigskip\begin{center}
\hspace{-.50cm}
\begin{overpic}[width=.46\textwidth]{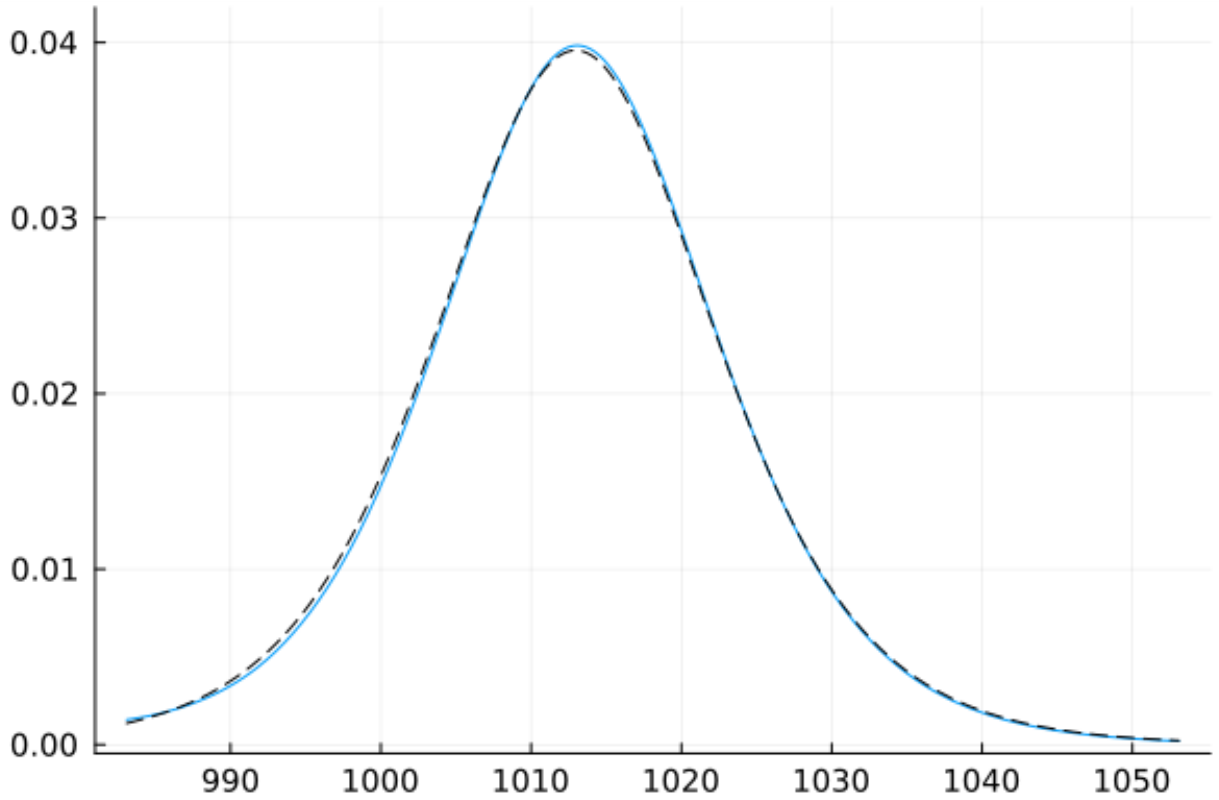}
      \put(0,68){\footnotesize $u_{\sol}, U$}
      \put(101.6,2.2){\footnotesize $x$}
       \put(40,25){\footnotesize $t = 1000$}
    \end{overpic}
        \hspace{0.4cm}
\begin{overpic}[width=.46\textwidth]{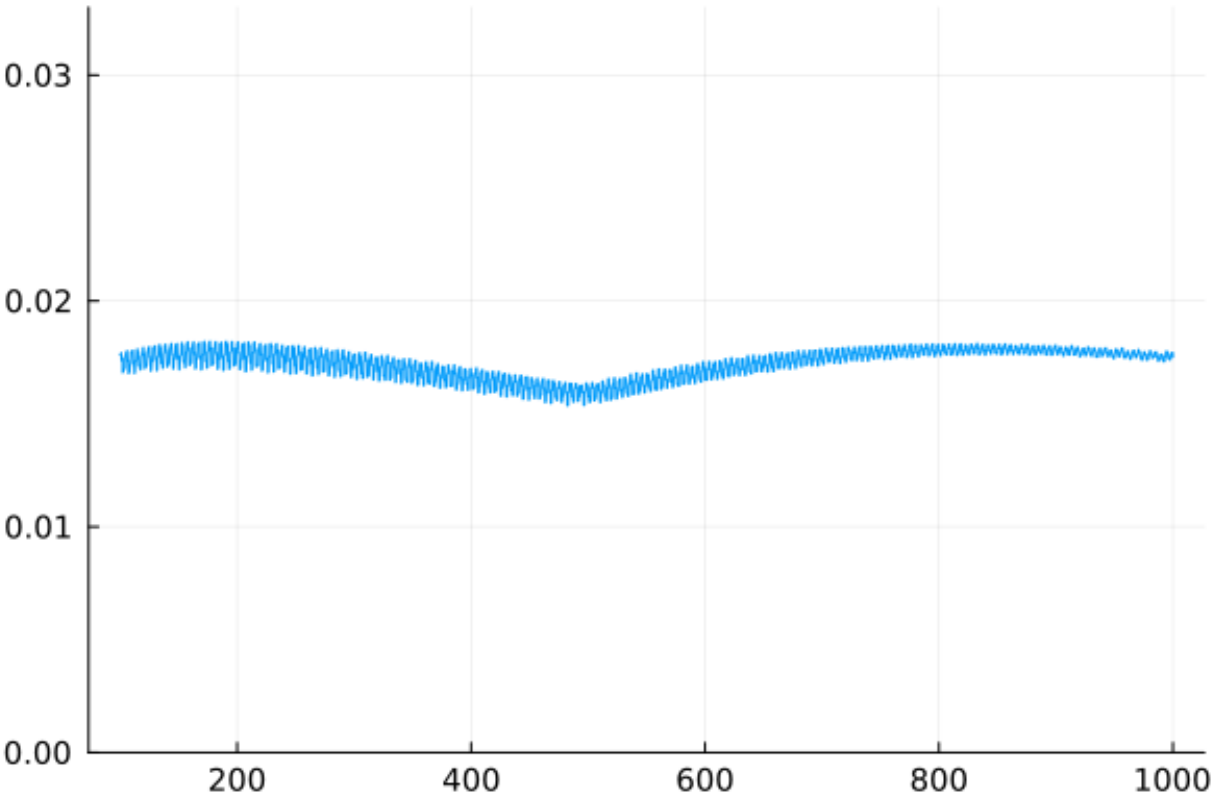}
      \put(0,68){\footnotesize $\sqrt{t}\,e(t)$}
      \put(101.6,2.2){\footnotesize $t$}
    \end{overpic}
     \begin{figuretext}\label{solutionwithsolitonfig3}
       Left: The asymptotic approximation $u_{\sol}(x,t)$ at $t = 1000$ (dashed black) together with the numerical approximation $U(x,t)$ (solid blue) for the initial data (\ref{solutionwithsolitoninitialdata}) obtained with the scheme of Section \ref{schemesec}. 
       
       Right: The graph of $t \mapsto \sqrt{t}\, e(t)$ where $e(t) = \sup_{x \in [1.001t, 2t]} |U(x,t) - u_{\sol}(x,t)|$.
\end{figuretext}
     \end{center}
\end{figure}

Figure \ref{solutionwithsolitonfig3} (left) shows $u_{\sol}(x,t)$ and $U(x,t)$ at $t = 1000$. Equation (\ref{uleadingsectorII}) with $\epsilon = 0.001$ and $M = 2$ implies that $u_{\sol}(x,t) - u(x,t)$ is of order $O(1/\sqrt{t})$ uniformly for $x/t \in [1.001t, 2t]$ as $t \to +\infty$. Figure \ref{solutionwithsolitonfig3} (right) shows the graph of the $L^\infty$-error $e(t) = \sup_{x \in [1.001t, 2t]} |U(x,t) - u_{\sol}(x,t)|$ multiplied by $\sqrt{t}$. If $U(x,t)$ is a good approximation of $u(x,t)$, then $\sqrt{t} e(t)$ should be of $O(1)$ at $t \to +\infty$. Figure \ref{solutionwithsolitonfig3} suggests that $\sqrt{t} e(t)$ is indeed of $O(1)$, confirming that $U(x,t)$ is a good approximation of $u(x,t)$.

\subsection{Conjecture}\label{conjecturesubsec}
Figure \ref{onehumpfig2}, Figure \ref{threehumpsfig2}, and Figure \ref{solutionwithsolitonfig2} suggest that the error term in (\ref{uasymptotics}) is uniformly valid for $\zeta$ near $1/\sqrt{3}$. This is noteworthy because analytically it has only been shown that (\ref{uasymptotics}) is valid uniformly for $\zeta$ in compact subsets of $(0, \frac{1}{\sqrt{3}}) \cup (\frac{1}{\sqrt{3}},1)$, see \cite{CLmain, CLsectorV, CLsectorIV}. Based on Figures \ref{onehumpfig2}, \ref{threehumpsfig2}, and \ref{solutionwithsolitonfig2}, as well as other numerical simulations, we conjecture that (\ref{uasymptotics}) holds uniformly for $\zeta$ in compact subsets of $(0, 1)$. In other words, we conjecture that there is no transition sector separating Sector IV and Sector V.

\appendix 

\section{The asymptotic formula (\ref{uleadingsectorII})}
We explain how the asymptotic formula (\ref{uleadingsectorII}) follows from the results of \cite{CLscatteringsolitons, CLsectorII}.

Define $\mathsf{U}(x,k)$, $X(x,k)$, $\mathcal{L}(k)$, $s(k)$, and $r_1(k)$ as in Section \ref{comparisonsubsubsec} with $u_0(x) = u(x,0)$ and $u_1(x) = u_t(x,0)$ given by (\ref{solutionwithsolitoninitialdata}).
Then $s_{11}(k)$ has a simple zero at $k_0 \approx 1.1755 >1$ and it was proved in \cite{CLscatteringsolitons} that there exists a unique complex constant $c_{k_0}$ such that 
$$\frac{[Y(x,k_0)]_2}{\dot{s}_{22}^A(k_0)} = c_{k_0} e^{(l_1(k_0)-l_2(k_0))x} [X(x,k_0)]_1\qquad \text{for all $x \in \R$},$$
where $Y(x,k)$ is the unique solution of the linear integral equation
\begin{align*}
 & Y(x,k) = I + \int_{-\infty}^x e^{(x-x')\mathcal{L}(k)} (\mathsf{U}Y)(x',k) e^{-(x-x')\mathcal{L}(k)} dx'.
\end{align*}
For $\zeta > 1$ and $k\in \C$ with $|k| \neq 1$, define
\begin{equation}
\delta(\zeta, k) = \exp \left\{ \frac{-1}{2\pi i} \int_{i}^{k_{1}} \frac{\ln(1 + \tilde{r}(s) |r_1(s)|^2)}{s - k} ds \right\},\qquad \tilde{r}(s) := \frac{\omega^2 - s^2}{1 - \omega^2 s^2}, 
\end{equation}
where the principal branch is taken for the logarithm and the path of integration follows the unit circle in the counterclockwise direction. Set
\begin{align*}
\hat{\Delta}_{33}(\zeta,k) := \frac{\delta(\zeta,\omega k)}{\delta(\zeta,\omega^{2} k)}\frac{\delta(\zeta,\frac{1}{\omega^{2} k})}{\delta(\zeta,\frac{1}{\omega k})}
\end{align*}
and consider the following Riemann--Hilbert (RH) problem. 

\begin{RHproblem}[RH problem for $M_{\mathrm{sol}}$]\label{RHS}
Find a $3 \times 3$-matrix valued function $M_{\mathrm{sol}}(x,t,k)$ with the following properties:
\begin{enumerate}[$(a)$]
\item\label{RHSitema} $M_{\mathrm{sol}}(x,t,\cdot) : \C \setminus \hat{\mathsf{Z}} \to \mathbb{C}^{3 \times 3}$ is analytic, where $\hat{\mathsf{Z}} = \{k_0,\omega k_0, \omega^2 k_0, k_0^{-1},\omega k_0^{-1}, \omega^2 k_0^{-1}\}$.

\item\label{RHSitemb} For $k \in \C \setminus  \hat{\mathsf{Z}}$, $M_{\mathrm{sol}}$ obeys the symmetries
\begin{align}\label{Msolsymm}
M_{\mathrm{sol}}(x,t,k) = \mathcal{A} M_{\mathrm{sol}}(x,t,\omega k)\mathcal{A}^{-1} = \mathcal{B} M_{\mathrm{sol}}(x,t,k^{-1}) \mathcal{B},
\end{align}
where
\begin{align}\label{def of Acal and Bcal}
\mathcal{A} := \begin{pmatrix}
0 & 0 & 1 \\
1 & 0 & 0 \\
0 & 1 & 0
\end{pmatrix} \qquad \mbox{ and } \qquad \mathcal{B} := \begin{pmatrix}
0 & 1 & 0 \\
1 & 0 & 0 \\
0 & 0 & 1
\end{pmatrix}.
\end{align}

\item\label{RHSitemc} $M_{\mathrm{sol}}(x,t,k) = I + O(k^{-1})$ as $k \to \infty$.

\item\label{RHSitemd}
At the zero $k_0 >1$ of $s_{11}(k)$, the first and third columns of $M_{\mathrm{sol}}$ are analytic, while the second column has (at most) a simple pole. Furthermore,
\begin{align}\label{Mk0inR}
 \underset{k = k_{0}}{\res} [M_{\mathrm{sol}}(x,t,k)]_2 = \frac{c_{k_{0}} e^{-\theta_{21}(x,t,k_0)}}{\hat{\Delta}_{33}(\zeta, \omega k_0)\hat{\Delta}_{33}^{-1}(\zeta, \omega^2 k_0)} [M_{\mathrm{sol}}(x,t,k_{0})]_1,
\end{align}
where $\theta_{21}(x,t,k)=(l_2(k)-l_1(k))x +(z_2(k)-z_1(k)) t$ with $z_{j}(k) = i \frac{(\omega^{j}k)^{2} + (\omega^{j}k)^{-2}}{4\sqrt{3}}$.
\end{enumerate}
\end{RHproblem}

It was proved in \cite[Theorem 2.1]{CLsectorII} that
\begin{align}\label{uusolurad}
& u(x,t) = u_{\mathrm{sol}}(x,t) + \frac{u_{\mathrm{rad}}(x,t)}{\sqrt{t}} + O\bigg(\frac{\ln t}{t}\bigg) \qquad \text{as $t\to +\infty$},
\end{align}
uniformly for $\zeta = \frac{x}{t}$ in compact subsets of $(1, +\infty)$, where the leading term $u_{\mathrm{sol}}$ is given in terms of the unique solution $M_{\mathrm{sol}}$ of RH problem \ref{RHS} by 
\begin{align}\label{usolMsol}
u_{\mathrm{sol}}(x,t) := -i\sqrt{3}\frac{\partial}{\partial x} \lim_{k\to\infty}k\Big(\big((1,1,1)M_{\mathrm{sol}}(x,t,k)\big)_{3}-1\Big).
\end{align}
The function $u_{\mathrm{rad}}(x,t)$ in (\ref{uusolurad}) is given explicitly in \cite[Theorem 2.1]{CLsectorII}, but for our present purposes it suffices to note that \cite[Lemmas A.1 and A.2]{CLsectorII} imply that $u_{\mathrm{rad}}(x,t)$ is of order $O(1)$ as $t \to +\infty$. To complete the derivation of (\ref{uleadingsectorII}), it therefore suffices to show that the function $u_{\mathrm{sol}}(x,t)$ defined in (\ref{usolMsol}) can be expressed as in (\ref{usoldef}). 

The fact that $u_{\mathrm{sol}}(x,t)$ can be expressed as in (\ref{usoldef}) follows from the exact solution of the above Riemann--Hilbert problem presented in \cite[Appendix A.1]{CLscatteringsolitons}. Indeed, applying the result of \cite[Appendix A.1]{CLscatteringsolitons} with $c_{k_0}$ replaced by $\frac{c_{k_{0}}}{\hat{\Delta}_{33}(\zeta, \omega k_0)\hat{\Delta}_{33}^{-1}(\zeta, \omega^2 k_0)}$, we find that
\begin{align}\label{uonezerodef}
u_\sol(x,t) = \frac{\frac{3}{2}(k_0 - k_0^{-1})^2}{\Big(f_{k_0}(\zeta) e^{-\frac{k_0^2 - 1}{4k_0}(x - c_{0}t)} + f_{k_0}(\zeta)^{-1} e^{\frac{k_0^2 - 1}{4k_0}(x - c_{0}t)}\Big)^2},
\end{align}
where
\begin{align}\label{cfk0def}
c_0 := \frac{k_0 + k_0^{-1}}{2}, \qquad  f_{k_0}(\zeta) := \sqrt{\frac{i \omega^2  (k_0^2 - \omega^2) c_{k_{0}}}{\sqrt{3}k_0(k_0^2 -1)} \frac{1}{\hat{\Delta}_{33}(\zeta, \omega k_0)\hat{\Delta}_{33}^{-1}(\zeta, \omega^2 k_0)}},
\end{align}
It is not necessary to specify the branch of the square root in (\ref{cfk0def}) because the value of $u(x,t)$ does not depend on the sign of $f_{k_0}(\zeta)$. By \cite[Lemma A.3]{CLscatteringsolitons}, we have $i \omega^2  (k_0^2 - \omega^2) c_{k_{0}} \geq 0$. The symmetry $\hat{\Delta}_{33}(\zeta,k) = \overline{\hat{\Delta}_{33}^{-1}(\zeta,\bar{k})}$ which is a consequence of \cite[Eq. (6.7)]{CLsectorII}, then implies that $f_{k_0}(\zeta)^2 \in \R$. In fact, since $u_\sol(x,t)$ is non-singular and non-trivial, we have $f_{k_0}(\zeta)^2 > 0$, so we can define a real-valued function $\ln f_{k_0}(\zeta)$ by
$$\ln f_{k_0}(\zeta) := \frac{1}{2} \ln\bigg(\frac{i \omega^2  (k_0^2 - \omega^2) c_{k_{0}}}{\sqrt{3}k_0(k_0^2 -1)} \frac{1}{\hat{\Delta}_{33}(\zeta, \omega k_0)\hat{\Delta}_{33}^{-1}(\zeta, \omega^2 k_0)}\bigg).$$
With this choice of $\ln f_{k_0}(\zeta)$, $u_{\mathrm{sol}}(x,t)$ can be expressed as in (\ref{usoldef}).

\subsection*{Acknowledgements}
Support is acknowledged from the Swedish Research Council, Grant No. 2021-04626 and Grant No. 2021-03877.

\bibliographystyle{plain}
\bibliography{is}

\begin{thebibliography}{99}
\small

\bibitem{BFH2019}
R. F. Barostichi, R. O. Figueira, and A. A. Himonas, Well-posedness of the ``good'' Boussinesq equation in analytic Gevrey spaces and time regularity, {\it J. Diff. Eq.} {\bf 267} (2019), 3181--3198. 

\bibitem{BS1988}
J. L. Bona and R. L. Sachs, Global existence of smooth solutions and stability of solitary waves for a generalized Boussinesq equation, {\it Comm. Math. Phys.} {\bf 118} (1988), 15--29. 

\bibitem{B1872}
J. Boussinesq, Th\'eorie des ondes et des remous qui se propagent le long d'un canal rectangulaire horizontal, en communiquant au liquide contenu dans ce canal des vitesses sensiblement pareilles de la surface au fond, {\it J. Math. Pures Appl.} {\bf 17} (1872), 55--108. 

\bibitem{B1998}
A. G. Bratsos, The solution of the Boussinesq equation using the method of lines, {\it Comput. Methods Appl. Mech. Engrg.}  {\bf 157} (1998), 33--44.

\bibitem{B2007}
A. G. Bratsos, A second order numerical scheme for the solution of the one-dimensional Boussinesq equation, {\it Numer. Algorithms}  {\bf 46} (2007), 45--58.

\bibitem{BTN2005}
A. G. Bratsos, Ch. Tsitouras, and D. G. Natsis, Linearized numerical schemes for the Boussinesq equation, {\it Appl. Numer. Anal. Comput. Math.}  {\bf 2} (2005), 34--53.

\bibitem{Cbad} C. Charlier, Blow-up solutions of the ``bad" Boussinesq equation, arXiv:2405.12210.

\bibitem{CL2022}
C. Charlier and J. Lenells, The ``good'' Boussinesq equation: a Riemann-Hilbert approach, {\it Indiana Univ. Math. J.}  {\bf 71} (2022), 1505--1562.

\bibitem{CLmain} 
C. Charlier and J. Lenells, On Boussinesq's equation for water waves, arXiv:2204.02365.

\bibitem{CLsectorV} 
C. Charlier and J. Lenells, Boussinesq's equation for water waves: asymptotics in Sector V, arXiv:2301.10669 (to appear in \textit{SIAM J. Math. Anal.}).

\bibitem{CLscatteringsolitons} 
C. Charlier and J. Lenells, Direct and inverse scattering for the Boussinesq equation with solitons, arXiv:2302.14593.

\bibitem{CLsectorIV} 
C. Charlier and J. Lenells, Boussinesq's equation for water waves: the soliton resolution conjecture for Sector IV, arXiv:2303.00434.

\bibitem{CLsectorII}
C. Charlier and J. Lenells, The soliton resolution conjecture for the Boussinesq equation, arXiv:2303.10485, 43pp.

\bibitem{CLW2022}
C. Charlier, J. Lenells, and D. Wang, The ``good" Boussinesq equation: long-time asymptotics, \textit{Anal. PDE \textbf{16} (2023), no.6, 1351--1388.}

\bibitem{CT2017}
E. Compaan and N. Tzirakis, Well-posedness and nonlinear smoothing for the ``good'' Boussinesq equation on the half-line, {\it J. Differential Equations} {\bf 262} (2017), 5824--5859. 

\bibitem{DH1999}
P. Daripa and W. Hua, A numerical study of an ill-posed Boussinesq equation arising in water waves and nonlinear lattices: filtering and regularization techniques, {\it Appl. Math. Comput.} {\bf 101} (1999), 159--207. 

\bibitem{D2008}
P. Deift, Some open problems in random matrix theory and the theory of integrable systems. In {\it Integrable systems and random matrices}, 419--430, Contemp. Math. {\bf 458}, Amer. Math. Soc., Providence, RI, 2008.

\bibitem{E2003}
H. El-Zoheiry, Numerical investigation for the solitary waves interaction of the ``good'' Boussinesq equation, {\it Appl. Numer. Math.} {\bf 45} (2003), 161--173.

\bibitem{FST1983}
G. E. Fal'kovich, M. D. Spector, and S. K. Turitsyn, Destruction of stationary solutions and collapse in the nonlinear string equation, {\it Phys. Lett. A} {\bf 99} (1983), 271--274.

\bibitem{F2009}
L. G. Farah, Local solutions in Sobolev spaces with negative indices for the ``good'' Boussinesq equation, {\it Comm. Partial Differential Equations} {\bf 34} (2009), 52--73. 

\bibitem{F2011}
L. C. F. Ferreira, Existence and scattering theory for Boussinesq type equations with singular data, {\it J. Diff. Eq.} {\bf 250} (2011), 2372--2388. 
 
\bibitem{HM2015}
A. A. Himonas and D. Mantzavinos, The ``good'' Boussinesq equation on the half-line, {\it J. Differential Equations} {\bf 258} (2015), 3107--3160.

\bibitem{H1973}
R. Hirota, Exact $N$-soliton solutions of the wave equation of long waves in shallow-water and in nonlinear lattices, {\it J. Math. Phys.} {\bf 14} (1973), 810--814. 

\bibitem{IB2003}
M. S. Ismail and A. G. Bratsos, A predictor-corrector scheme for the numerical solution of the Boussinesq equation, {\it J. Appl. Math. Comput.} {\bf 13} (2003), 11--27.

\bibitem{IM2014}
M. S. Ismail and F. Mosally, A fourth order finite difference method for the good Boussinesq equation, {\it Abstr. Appl. Anal.} (2014), Art. ID 323260, 10 pp.

\bibitem{J1997}
R. S. Johnson, A modern introduction to the mathematical theory of water waves. Cambridge Texts in Applied Mathematics. Cambridge University Press, Cambridge, 1997.

\bibitem{KL1977}
V. K. Kalantarov and O. A. Ladyzenskaja, Formation of collapses in quasilinear equations of parabolic and hyperbolic types. (Russian) Boundary value problems of mathematical physics and related questions in the theory of functions, 10. {\it Zap. Nau\v{c}n. Sem. Leningrad. Otdel. Mat. Inst. Steklov (LOMI)} {\bf 69} (1977), 77--102.

\bibitem{KT2010} 
N. Kishimoto and K. Tsugawa, Local well-posedness for quadratic nonlinear Schr\"odinger equations and the ``good'' Boussinesq equation, {\it 
Differential Integral Equations} {\bf 23} (2010), 463--493. 

\bibitem{L1993}
F. Linares, Global existence of small solutions for a generalized Boussinesq equation, {\it J. Differential Equations} {\bf 106} (1993), 257--293. 

\bibitem{LS1995}
F. Linares and M. Scialom, Asymptotic behavior of solutions of a generalized Boussinesq type equation, {\it Nonlinear Anal.} {\bf 25} (1995), 1147--1158.

\bibitem{L1997}
Y. Liu, Decay and scattering of small solutions of a generalized Boussinesq equation, J. Funct. Anal. {\bf 147} (1997), 51--68. 

\bibitem{M1978}
V. G. Makhankov, Dynamics of classical solitons (in nonintegrable systems), {\it Phys. Rep.} {\bf 35} (1978), 1--128.
 
\bibitem{MMM1985}
V. S. Manoranjan, A. R. Mitchell, and J. Ll. Morris, Numerical solution of the good Boussinesq equation, {\it SIAM J. Sci. Stat. Comput.} {\bf 5} (1984), 946--957.
 
  \bibitem{OT2013} 
 S. Olver, and A. Townsend, A fast and well-conditioned spectral method, \textit{SIAM Rev.} \textbf{55} (2013), 462--489.

\bibitem{OS1990}
T. Ortega and J. M. Sanz-Serna, Nonlinear stability and convergence of finite-difference methods for the ``good'' Boussinesq equation, {\it Numer. Math.} {\bf 58} (1990), 215--229.

\bibitem{PS1997}
A. K. Pani and H. Saranga, Finite element Galerkin method for the ``good'' Boussinesq equation, {\it Nonlinear Anal.} {\bf 29} (1997), 937--956.

\bibitem{T1975}
M. Toda, Studies of a non-linear lattice, {\it Phys. Rep.} {\bf 18C} (1975), 1--123. 

\bibitem{TM1991}
M. Tsutsumi and T. Matahashi, On the Cauchy problem for the Boussinesq type equation, {\it Math. Japon.} {\bf 36} (1991), 371--379.

\bibitem{UEK2021}
Y. Ucar, A. Esen, and B. Karaagac, Numerical solutions of Boussinesq equation using Galerkin finite element method, {\it Numer. Methods Partial Differential Equations} {\bf 37} (2021), 1612--1630.

\bibitem{X2006}
R. Xue, Local and global existence of solutions for the Cauchy problem of a generalized Boussinesq equation, {\it J. Math. Anal. Appl.} {\bf 316} (2006), 307--327. 

\bibitem{Z1974}
V. E. Zakharov, On stochastization of one-dimensional chains of nonlinear oscillations, {\it Soviet Phys. JETP} {\bf 38} (1974), 108--110.

\end{thebibliography}

\end{document}